\documentclass[11pt,a4paper]{article}
\usepackage{amsfonts}
\usepackage{amsmath}
\usepackage{amssymb}
\usepackage{amsthm}
\usepackage{enumitem}
\usepackage{hyperref}
\usepackage{mathtools}
\usepackage{graphicx}
\usepackage{color}
\usepackage[left=2.9cm,right=2.9cm, bottom=3.5cm, top=3.5cm]{geometry}
\usepackage[all]{xy}
\usepackage{tikz}
\usetikzlibrary{decorations.markings}
\usetikzlibrary{decorations.pathreplacing}
\usetikzlibrary{arrows,shapes,positioning}
\tikzstyle directed=[postaction={decorate,decoration={markings,    mark=at position #1 with {\arrow{>}}}}]
\tikzstyle rdirected=[postaction={decorate,decoration={markings,   mark=at position #1 with {\arrow{<}}}}]
\allowdisplaybreaks
\tikzset{anchorbase/.style={baseline={([yshift=-0.5ex]current bounding box.center)}}}
\tikzset{anchorbase/.style={baseline={([yshift=-0.5ex]current bounding box.center)}}}

\newcommand\scalemath[2]{\scalebox{#1}{\mbox{\ensuremath{\displaystyle #2}}}}

\usepackage{bbm}
\def\C{{\mathbbm C}}
\def\N{{\mathbbm N}}
\def\R{{\mathbbm R}}
\def\Z{{\mathbbm Z}}

\def\cal#1{\mathcal{#1}}%

\def\comm#1{}%

\newcommand{\qp}[2]{(#1;q^2)_{#2}}
\newcommand{\qpp}[1]{(q^2;q^2)_{#1}}

\theoremstyle{definition}
\newtheorem{theorem}{Theorem}[section]
\newtheorem{corollary}[theorem]{Corollary}

\newtheorem{lemma}[theorem]{Lemma}
\newtheorem{remark}[theorem]{Remark}
\newtheorem{proposition}[theorem]{Proposition}
\newtheorem{definition}[theorem]{Definition}
\newtheorem{example}[theorem]{Example}

\newtheorem{notation}[theorem]{Notation}

\title{Rational links and DT invariants of quivers}

\author{Marko Sto$\check{\text{s}}$i$\acute{\text{c}}$
  \and
  Paul Wedrich
 }

\date{}

\begin{document}

\maketitle
\abstract{We prove that the generating functions for the colored HOMFLY-PT polynomials of rational links are specializations of the generating functions of the motivic Donaldson-Thomas invariants of appropriate quivers that we naturally associate with these links. This shows that the conjectural links-quivers correspondence of Kucharski--Reineke--Sto$\check{\text{s}}$i$\acute{\text{c}}$--Su{\l}kowski as well as the LMOV conjecture hold for rational links. Along the way, we extend the links-quivers correspondence to tangles and, thus, explore elements of a skein theory for motivic Donaldson-Thomas invariants.}

\section{Introduction}

Consider a symmetric quiver with $n$ vertices in a chosen order and adjacency matrix $Q$. Efimov~\cite{E} has shown that the cohomological Hall algebra associated to this symmetric quiver (without potential) by Kontsevich--Soibelman~\cite{KS} is free-supercommutative (after twisting the multiplication by a sign). Its generating function, the $\Z\times \N^n$-graded Hilbert--Poincar\'{e} series, is given by:

\begin{equation}\label{eqn:motivic}
P_Q(\textbf{x})= \sum_{\textbf{d}\in \N^n} (-q)^{-\langle \textbf{d},\textbf{d} \rangle_Q} \textbf{x}^{\textbf{d}} \prod_{i=1}^n \prod_{k=1}^{d_i} \frac{1}{1-q^{-2k}} \in \Z(q)[[x_1,\ldots,x_n]]
\end{equation} 
Here $\textbf{d}=(d_1,\dots,d_n)$, $\textbf{x}^{\textbf{d}}=x_1^{d_1}\cdots x_n^{d_n}$ and $\langle \textbf{d},\textbf{e} \rangle_Q:=\sum_{i,j}(\delta_{i,j}-Q_{i,j})d_i e_j$ is the Euler form of the quiver.
\newline 

The reduced $j$-colored HOMFLY-PT polynomial $P_j$ is an invariant of framed, oriented links in $\R^3$ that takes values in $\Z[a^{\pm 1}](q)$. It can be computed as a limit of Reshetikhin--Turaev/Chern--Simons invariants in type A at large rank and the normalization is chosen such that $P_j(\bigcirc)=1$. The related Labastida--Mari\~no--Ooguri--Vafa (LMOV) or BPS invariants of links are a repackaging of colored HOMFLY-PT invariants motivated by topological string theory, see \cite{OV,LM}.\\

Recently, a surprising relationship between 
 invariants of links and motivic Donaldson--Thomas quiver invariants has been
 conjectured in \cite{KRSSlong, KRSSshort}. More precisely, the conjecture
 posits that for each link there exists a quiver, such that its generating
 function \eqref{eqn:motivic} determines the generating function of the
 $j$-colored HOMFLY-PT polynomials of the link. In \cite{KRSSlong} such quivers
 have been obtained for all knots with at most six crossings, as well as all
 $(2,2n+1)$ torus knots and twist knots. We refer the interested reader to these
 sources for more details about the motivation and background for this
 relationship.\\

The main purpose of the present paper is to prove the links-quivers correspondence conjectured in \cite{KRSSlong,KRSSshort} for all rational (a.k.a. 2-bridge) links.

\begin{theorem}\label{thm:knots} Let $K_{u/v}$ be a positive rational knot. Then there exists a quiver $Q_{K_{u/v}}$ with $u$ vertices and $A=(a_1,\dots,a_u),  S=(q_1,\dots,q_u)\in \Z^u$ such that:
\[\sum_{j\geq 0} \frac{P_j(K_{u/v})}{\prod_{i=1}^{j}(1-q^{2i})} x^j = P_{Q_{K_{u/v}}}(\textbf{x})\mid_{x_i \mapsto (-1)^{Q_{i,i}+q_i} q^{q_i-1}a^{a_i} x} \]
Recall that $P_j$ are normalized (reduced) $j$-colored HOMFLY-PT polynomials of the knot $K_{u/v}$, and $P_{Q_{K_{u/v}}}(\textbf{x})$ is the generating function (\ref{eqn:motivic}) of the quiver $Q_{K_{u/v}}$.
\end{theorem}

For 2-component links, the optimal analogous result features modified generating functions, see Remark~\ref{rem:HOMFLY-links}.  
\begin{theorem}
\label{thm:links} Let $L_{u/v}$ be a positive rational (2-component) link, then there exists a quiver $Q_{L_{u/v}}$ with $2u$ vertices and integer vectors $A=(a_1,\dots,a_{2u})$ and $S=(q_1,\dots,q_{2u})$, such that the generating function of the reduced $j$-colored HOMFLY-PT polynomials can be obtained from the generating function of $Q_{L_{u/v}}$ by specializing $x_i \mapsto (-1)^{Q_{i,i}+q_i} q^{q_i-1}a^{a_i} x$:
\[\sum_{j\geq 0} P_j(L_{u/v}) x^j = P_{Q_{L_{u/v}}}(\textbf{x})\mid_{x_i \mapsto (-1)^{Q_{i,i}+q_i} q^{q_i-1}a^{a_i} x} \]
\end{theorem}

One of the main applications of the links-quivers correspondence is that the LMOV invariants of a link can be written as an integral linear combination of the quantum DT invariants of the corresponding quiver. Since these are integral \cite[Section 6.2]{KS} (and indeed positive \cite[Section 4]{E}, with important consequences for the positivity of Kac polynomials \cite{Moz}, see also \cite{HLV,D}), this implies the integrality of the LMOV invariants of the corresponding link, see \cite{KRSSlong}.
\begin{corollary}
The LMOV integrality conjecture holds for all rational links with (anti-)symmetric colors.
\end{corollary}

Another observation of \cite{KRSSlong} is that the quivers $Q_K$ for knots $K$ often encode information about the reduced (colored) HOMFLY-PT {\it{homology}} of $K$, \cite{KR, Wed3}. We confirm their main conjecture in the case of rational knots. 

\begin{theorem}\label{thm:knotsHHH} Let $K$ be a rational knot and let $Q_K$ be the corresponding quiver from Theorem \ref{thm:knots}. Then, the vertices of $Q_K$ are in bijection with generators of the reduced HOMFLY-PT homology of $K$, such that the $(a,q,t)$-trigrading of the $i^{th}$ generator is given by $(a_i,-Q_{i,i}-q_i,-Q_{i,i})$ where $Q_{i,i}$ denotes the number of loops at the $i^{th}$ vertex of $Q_K$.
\end{theorem}

In order to prove Theorems \ref{thm:knots}, \ref{thm:links} and \ref{thm:knotsHHH}, we use an algorithm from \cite{Wed2} for the computation of the colored HOMFLY-PT polynomials of rational links. The algorithm iteratively computes the HOMFLY-PT invariants of rational tangles, one crossing at a time, before performing a closure operation to arrive at the invariant of the rational link. A key observation in the present paper is that not only rational links, but also rational \textit{tangles} have associated quivers whose generating functions describe their HOMFLY-PT invariants. Moreover, the elementary operations of adding crossings to rational tangles correspond to natural operations on quivers, which double certain subsets of vertices. Using these operations, the quivers for rational links can be computed iteratively from a continued fraction expansion. It is an intriguing problem for further research to extend this skein theory for motivic DT invariants beyond the realm of rational tangles and links.

\begin{remark} Relationships between HOMFLY-PT and motivic DT invariants have been considered by Diaconescu--Hua--Soibelman~\cite{DHS} in the context of the Oblomkov--Shende conjecture~\cite{OS}, which was proved by Maulik~\cite{Mau}. For an extension to HOMFLY-PT homology, see Oblomkov--Rasmussen--Shende~\cite{ORS}. These results concern algebraic links, i.e. links of plane curve singularities, which intersect the class of rational links only in the 2-strand torus links. The problem of constructing quivers corresponding to algebraic links beyond small torus links remains open.
\end{remark}

 \textbf{Organisation of the paper.} Section~\ref{sec:background} contains background information on rational links and tangles, the skein theory of the colored HOMFLY-PT invariants and a few useful formulas for manipulating generating functions. In Section~\ref{sec:links} we present an algorithm for computing quivers for rational links. In Section~\ref{sec:knots} we refine this algorithm to reduce the sizes of these quivers and prove Theorems~\ref{thm:knots} and~\ref{thm:links} as well as Theorem~\ref{thm:knotsHHH}. Section~\ref{sec:examples} contains example computations.\\

\section{Background}\label{sec:background}
We start by introducing notation and terminology and then recall background about rational tangles and their HOMFLY-PT skein theory.

\subsection{Quantum algebra}
We will use the Pochhammer symbols $(x;y)_i=\prod_{j=0}^{i-1}(1-x y^j)$ and especially the $q$-Pochhammer symbols $\qpp{i}=\prod_{j=1}^{i}(1-q^{2j})$. For non-negative integers $a_1,\ldots,a_k$ with $a_1+\ldots+a_k=N$, we define the quantum multinomial coefficients:
\begin{equation}
{
N\brack
a_1,\ldots,a_k
}=
\frac{(q^2;q^2)_N}{(q^2;q^2)_{a_1}\cdots(q^2;q^2)_{a_k}}.
\end{equation}
These are polynomials in $q^2$ with constant term 1. In the special case $k=2$ they are the (positive) quantum binomial coefficient:
$${N\brack k,N-k}={N \brack k}_+$$
and we will also use this notation.

\begin{definition}
\label{def:linky} We say that a generating function $P(x)\in \Z[a^{\pm 1}](q)[[x]]$ is in \textit{quiver form} if it can be obtained from a generating function \eqref{eqn:motivic} by specializing the variables $x_i \mapsto (-1)^{Q_{i,i}+q_i} q^{q_i-1}a^{a_i} x$ for some integer vectors $A=(a_1,\dots,a_n)$ and $S=(q_1,\dots,q_n)$ and an $n\times n$ symmetric matrix $Q$ with entries in $\Z$. The resulting expression is of the form
\begin{equation}
\label{eqn:quiverform}
P(x)=\sum_{\textbf{d}=(d_1,\dots,d_n)\in \N^n} \frac{(-q)^{S\cdot\textbf{d}} a^{A\cdot\textbf{d}} q^{\textbf{d}\cdot Q\cdot\textbf{d}^t} {x^{d_1+\cdots+d_n}} }{\qpp{d_1}\cdots \qpp{d_n}}  
\end{equation}
\end{definition}

If the generating function has an additional factor of $\qpp{d_1+\cdots+d_n}$ in each summand, the $q$-Pochhammer symbols combine to give a $q$-multinomial coefficient.
\begin{definition}
\label{def:knotty}
A generating function of the form

\begin{equation}
\label{eqn:knottyquiverform}
P(x)=\sum_{\textbf{d}=(d_1,\dots,d_n)\in \N^n} (-q)^{S\cdot\textbf{d}} a^{A\cdot\textbf{d}} q^{\textbf{d}\cdot Q\cdot\textbf{d}^t} x^{d_1+\cdots+d_n} {d_1+\cdots+d_n \brack d_1,\ldots,d_n},  
\end{equation}
for $S, A\in \Z^n$ and an $n\times n$ symmetric integer matrix $Q$, is said to be in \textit{polynomial quiver form}. 
\end{definition}
Using this terminology, Theorem~\ref{thm:knots} claims that the generating functions of colored HOMFLY-PT polynomials of rational knots can be brought into polynomial quiver form.

We say that a generating function is \textit{almost in quiver form}, if it is of the form \eqref{eqn:quiverform}, but potentially with several additional $q$-Pochhammer symbols
$\qpp{d_{i_1}+\cdots+d_{i_k}}$ as factors in the summands. Any generating function that is almost in quiver form can be re-written in quiver form by the following lemma.

\begin{lemma}[{\cite[Lemma 4.5]{KRSSlong}}]\label{lemlong}
For any $d_1,\ldots,d_k\ge 0$, we have:
\begin{align*}\frac{(x^2;q^2)_{d_1+\ldots+d_k}}{(q^2;q^2)_{d_1}\cdots (q^2;q^2)_{d_k}}
=& \!\!\!  \sum\limits_{\substack{\alpha_1+\beta_1=d_1\\\cdots\\\alpha_k+\beta_k=d_k}}\!\!\!\!\!\!\!\!  
\frac{(-x^2 q^{-1})^{{\alpha_1+\ldots+\alpha_k}} q^{\alpha_1^2+\ldots+\alpha_k^2+2\sum_{i=1}^{k-1} \alpha_{i+1} (d_1+\ldots+d_i)}}{(q^2;q^2)_{\alpha_1}\cdots(q^2;q^2)_{\alpha_k}(q^2;q^2)_{\beta_1}\cdots(q^2;q^2)_{\beta_k}}      
\end{align*}
\end{lemma}

This follows from the \textit{$q$-binomial identity:}
\begin{lemma}\label{mainlemma} We have
$$\qp{x^2}{k}=\sum_{i=0}^k (-1)^i x^{2i} q^{i^2 -i} {k \brack i}_+ = \sum_{i=0}^k (-1)^i x^{2i} q^{i^2 -i} \frac{\qpp{k}}{\qpp{i}\qpp{k-i}}.$$
\end{lemma}

Another useful tool for manipulating generating series is the following.

\begin{lemma}[{\cite[Lemma 4.6]{KRSSlong}}]\label{lemvlong}
Let $a_1,\ldots,a_m$ and $b_1,\ldots,b_p$ be nonnegative integers satisfying:
$$a_1+\ldots+a_m=b_1+\ldots+b_p.$$
Then we have:
\begin{equation}\label{forexp}
\!\!\!{a_1+\ldots+a_m \brack b_1,\ldots,b_p}=
\sum_{
\begin{array}{c}
\{j_{\alpha,\beta}\}\\ 
{\scriptstyle{{\alpha=1,\ldots,p}}} \\
{\scriptstyle{\beta=1,\ldots,m}}
\end{array}
}
q^{X(\underline{j})}{a_1 \brack j_{1,1},j_{2,1},\ldots,j_{p,1}}
\cdots {a_m \brack j_{1,m},j_{2,m},\ldots,j_{p,m}},
\end{equation}

\noindent with 
\begin{equation}
\label{eqn:triangcorr} X(\underline{j})=2\sum_{\scriptstyle{1\le l_1<l_2\le p}}\,\,\sum_{\scriptstyle{1\le u_1<u_2\le m}} \,j_{l_1,u_1} j_{l_2,u_2},
\end{equation}
where in formula (\ref{forexp}) we are summing over $m p$ summation indices (nonnegative integers) $j_{\alpha,\beta}$, with $\alpha=1,\ldots,p$ and $\beta=1,\ldots,m$, such that
\begin{eqnarray}
&j_{1,v}+ j_{2,v}+\ldots+ j_{p,v}=a_v,\quad v=1,\ldots,m,\\
&j_{w,1}+j_{w,2}+\ldots+j_{w,m}=b_w,\quad w=1,\ldots,p.
\end{eqnarray}
\end{lemma}

\subsection{Rational tangles and links}

Positive rational tangles are oriented 4-ended tangles, which can be built, starting from 
the trivial tangle $\;\begin{tikzpicture} [scale=.3,anchorbase]
\draw[thick, directed=.65] (-0.5,0) to (-0.5,1);
\draw[thick, directed=.65] (0.5,0) to (0.5,1);
\draw[thick, white] (0,-.1) to (.1,-.1);
\end{tikzpicture}\;$
by iteratively attaching a finite number of crossings $\begin{tikzpicture} [scale=.3,anchorbase]
\draw[thick] (0.5,0) to (-0.5,1);
\draw[white, line width=.15cm] (-0.5,0) to (0.5,1);
\draw[thick] (-0.5,0) to (0.5,1);
\draw[thick, white] (0,-.1) to (.1,-.1);
\end{tikzpicture}$
to the two top endpoints or the two right endpoints. We call these operations top and right twist respectively.

\[
T\left(\begin{tikzpicture} [scale=.5,anchorbase]
\draw[thick] (-0.5,0) to (-0.5,1.5);
\draw[thick] (0.5,0) to (0.5,1.5);
\draw[thick, fill=white] (-.75,.5) rectangle (.75,1);
\end{tikzpicture}\right)
:=
\begin{tikzpicture} [scale=.5,anchorbase]
\draw[thick] (0.5,0) to (0.5,1) to [out=90,in=270] (-.5,2);
\draw[white, line width=.15cm] (-0.5,1) to [out=90,in=270] (.5,2);
\draw[thick] (-0.5,0) to (-0.5,1) to [out=90,in=270] (.5,2);
\draw[thick, fill=white] (-.75,.5) rectangle (.75,1);
\end{tikzpicture}
\quad, \quad
R\left(\begin{tikzpicture} [scale=.5,anchorbase]
\draw[thick] (-0.5,0) to (-0.5,1.5);
\draw[thick] (0.5,0) to (0.5,1.5);
\draw[thick, fill=white] (-.75,.5) rectangle (.75,1);
\end{tikzpicture}\right)
:=
\begin{tikzpicture} [scale=.5,anchorbase]
\draw[thick] (-0.5,0) to (-0.5,1.5);
\draw[thick] (0.5,.5) to (0.5,1) to [out=90,in=180](0.7,1.2) to [out=0,in=135] (1.5,.75) to [out=325,in=90] (2,0);
\draw[white, line width=.15cm] (0.7,.3) to [out=0,in=225] (1.5,.75) to [out=45,in=270] (2,1.5);
\draw[thick] (0.5,.5) to [out=270,in=180] (0.7,.3) to [out=0,in=225] (1.5,.75) to [out=45,in=270] (2,1.5);
\draw[thick, fill=white] (-.75,.5) rectangle (.75,1);
\end{tikzpicture}
\]

We may assume that the building process starts with a positive number $a_1$ of top twists followed by a positive number $a_2$ of right twists. We continue to record the numbers $a_i$ of consecutive twists of the same type added to the tangle, where even indices $i$ encode right twists and odd indices top twists. The sequence of non-zero $a_i$ is then interpreted as the continued fraction expansion of a rational number. 

\[u/v = [a_1,\dots,a_r] = a_r + \frac{1}{a_{r-1}+\frac{1}{a_{r-2}+\dots}} \]

Note that for such non-trivial tangles $\tau_{u/v}$ we have $u/v>1$, and we assume $\mathrm{gcd}(u,v)=1$ as well as $u>0$ and $v>0$. We also record the length $r$ of the continued fraction expansion. As an example, we draw the tangle encoded by $[2,3,1]$:\
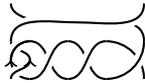
\[\begin{tikzpicture} [scale=.25,anchorbase]
\draw[thick] (1,1) to [out=90,in=270] (0,2);
\draw[thick] (0,1)to [out=270,in=180] (1,0);
\draw[thick] (7,0) to [out=90,in=0] (5,2);
\draw[thick] (3,0) to [out=180,in=0] (1,2);
\draw[thick] (3,2)to [out=0,in=180] (5,0);
\draw[thick] (7,2)to [out=90,in=0] (3.5,3)to [out=0,in=270] (0,4);
\draw[white, line width=.15cm] (0,2) to [out=90,in=180] (3.5,3) to[out=0,in=270] (7,4);
\draw[white, line width=.15cm] (0,0) to [out=90,in=270](1,1);
\draw[white, line width=.15cm] (1,2) to [out=180,in=90] (0,1);
\draw[white, line width=.15cm] (5,2)to [out=180,in=0] (3,0);
\draw[white, line width=.15cm] (1,2) to [out=180,in=90] (0,1);
\draw[white, line width=.15cm] (1,0) to [out=0,in=180] (3,2);
\draw[white, line width=.15cm] (5,0)to [out=0,in=270] (7,2);
\draw[thick] (0,2) to [out=90,in=180] (3.5,3) to[out=0,in=270] (7,4);
\draw[thick] (0,0) to [out=90,in=270](1,1);
\draw[thick] (1,2) to [out=180,in=90] (0,1);
\draw[thick] (5,2)to [out=180,in=0] (3,0);
\draw[thick] (1,2) to [out=180,in=90] (0,1);
\draw[thick] (1,0) to [out=0,in=180] (3,2);
\draw[thick] (5,0)to [out=0,in=270] (7,2);
\draw[thick,->](1,1) to(1,1.01);
\draw[thick,->](0,1) to(0,1.01);
\end{tikzpicture}\]

If the length $r$ is odd, i.e. the tangle building process ends with a top twist, we close the tangle using the NS-closure by connecting North to South boundary points in a planar fashion (provided the orientations match up). If the length $r$ is even, we close in the EW direction. The resulting knots and links are called positive rational or $2$-bridge, and they only depend on the fraction $u/v$, so we denote them by $L_{u/v}$. In particular, the two distinct continued fraction expansions $[a_1+1,\dots, a_r]=[1,a_1,\dots,a_r]$ of the same rational number give isotopic links as closures. Negative rational links can be obtained by taking mirror images.

The type of closure we have described is called the numerator closure, because the determinant of a link, given as the closure of a positive rational tangle with fraction $u/v$ is equal to $u$. The other possible way of closing, i.e. the one which leads to nugatory crossings, produces links with determinant $v$ and is known as the denominator closure.

The boundary orientations of a rational tangle can be in three different configurations, which we denote by $UP$, $OP$ and $RI$, which stands for upward, opposite and rightward.

\[UP:\; \begin{tikzpicture} [scale=.5,anchorbase]
\draw[thick, ->] (-0.5,0) to (-0.5,1.5);
\draw[thick, ->] (0.5,0) to (0.5,1.5);
\draw[thick, fill=white] (-.75,.5) rectangle (.75,1);
\end{tikzpicture}
\;,\quad 
OP: \;\begin{tikzpicture} [scale=.5,anchorbase]
\draw[thick, ->] (-0.5,0) to (-0.5,1.5);
\draw[thick, <-] (0.5,0) to (0.5,1.5);
\draw[thick, fill=white] (-.75,.5) rectangle (.75,1);
\end{tikzpicture}
\;,\quad
RI:\; \begin{tikzpicture} [scale=.5,anchorbase]
\draw[thick] (-0.5,0) to (-0.5,1.5);
\draw[thick, <->] (0.5,0) to (0.5,1.5);
\draw[thick, fill=white] (-.75,.5) rectangle (.75,1);
\end{tikzpicture}\]

 $UP$ tangles can only be closed in the NS direction, $RI$ tangles only in the EW direction, while $OP$ tangles admit both types of closures. 

Let $L_{u/v}$ and $L_{u'/v'}$ be two oriented rational links with $u>v>0$ and $u'>v'>0$. Then $L_{u/v} = L_{u'/v'}$ as unoriented links if and only if $u=u'$ and either $v=v'$ or $vv'\equiv 1 ~(\mathrm{mod }~ u)$. Furthermore, they are equal as oriented links if and only if $u=u'$ and either $v=v'$ or $vv'\equiv 1 ~(\mathrm{mod }~ 2u)$. We also write $L_{[a_1,\dots,a_r]}$ for $L_{u/v}$ if $[a_1,\dots,a_r]$ is a continued fraction expansion for $u/v$. Then $L_{[a_1,\dots,a_r]}$ and $L_{[a_r,\dots,a_1]}$ are isotopic as unoriented links and they are even oriented-isotopic in the case of knots. These facts are well-known, see e.g. \cite{Mur}.

\begin{lemma} The positive rational link $L_{u/v}$ is a knot if and only if $u$ is odd. It can be built via an odd-length continued fraction expansion 
ending with $UP$ if $v$ is odd and with $OP$ if $v$ is even.
\end{lemma} 

\begin{proof} 
We consider $u/v=[a_1,\dots,a_r]$ with odd $r$ and the link $L_{[a_1,\dots,a_r]}$. We now refine the notation $UP$, $OP$ and $RI$ for the boundary notation by including connectivity information between these boundaries: we let $UPk$ ($OPk$) and $UPl$ ($OPl$) denote $UP$ ($OP$) configurations which respectively give knots or 2-components links under NS-closure. Analogously, $RIk$ and $RIl$ denote $RI$ configurations that close to knots or 2-component links under EW-closure.
We associate $UPl$ and $OPl$ with the symbol $e/o$, $OPk$ and $RIl$ with $o/e$ and $UPk$ and $RIk$ with $o/o$. By assumption, the trivial starting tangle is of type $UPl$ and this corresponds to the rational number $0/1$, which is even/odd, as encoded by the symbol $e/o$. We can now check that applying top and right twists to these configurations leads to tangle types which correspond to partial continued fraction expansions of the prescribed parities. Note that individual top and right twists map between these configurations as follows:
\begin{align*}
 UPl
\xleftrightarrow{T} UPk  
\xleftrightarrow{R} OPk 
\xleftrightarrow{T} RIl  
\xleftrightarrow{R} RIk  
\xleftrightarrow{T} OPl  
\xleftrightarrow{R} UPl
\\
 e/o 
 \;\xleftrightarrow{}\; o/o 
 \;\xleftrightarrow{}\; o/e 
 \;\xleftrightarrow{}\; o/e 
 \;\xleftrightarrow{}\; o/o 
 \;\xleftrightarrow{}\; e/o 
 \;\xleftrightarrow{}\; e/o 
 \end{align*}
 After the first set of top twists, starting with $UPl$, we have the rational number $a_1/1$ and are in configuration $UPl$ or $UPk$ depending on whether $a_1$ is even or odd. This matches with the types $e/o$ and $o/o$ respectively. Next we apply $T^b R^a$ to the rational tangle associated to $x/y$, which produces the tangle encoded by: 
 \[b + \frac{1}{a+\frac{1}{x/y}} =(b a x + b y + x)/(a x +y)\] The first of the following two tables shows the parities of $(b a x + b y + x)/(a x +y)$ depending on the parities of $x/y$, $a$ and $b$, while the second table shows the orientation and connectivity of the tangle after applying $T^b R^a$, depending on the initial configuration.
 
 \[\scalemath{0.8}{\begin{tabular}{c c| c c c }
  b& a & e/o & o/e & o/o
 \\ \hline\hline
  e & e & e/o & o/e & o/o
   \\ \hline
  o & e & o/o & o/e & e/o
 \\ \hline
  e & o & e/o & o/o & o/e
   \\ \hline
  o & o & o/o & e/o & o/e
 \end{tabular}   
 \quad
 \begin{tabular}{c c| c c c }
  b& a & UPl,OPl & OPk,RIl & UPk,RIk
 \\ \hline\hline
  e & e & UPl,OPl & OPk,RIl & UPk,RIk
   \\ \hline
  o & e & UPk,RIk & RIl,OPk & UPl,OPl
 \\ \hline
  e & o & OPl,UPl & UPk,RIk & OPk,RIl
   \\ \hline
  o & o & RIk,UPk & UPl,OPl & RIl,OPk
 \end{tabular} 
 }  
\] 
Comparing these tables shows that the correspondence between the parities of the fraction and the orientation and connectivity of the tangle is preserved throughout the building process of the rational tangle. In particular, a NS-closure of the rational tangle encoded by the positive odd-length continued fraction expansion of $u/v$ gives a knot if and only if $u$ is odd, since these are exactly the cases when we can get $UPk$ and $OPk$ at the end. Moreover, $UPk$ precisely arises when $v$ is odd and $OPk$ when $v$ is even.
\end{proof}

\subsection{HOMFLY-PT skein theory}

We will consider the HOMFLY-PT invariants of rational links with all components colored by the one-column Young diagram with $j$ boxes. For a link $L$ we denote this unreduced $j$-colored HOMFLY-PT invariant by $\overline{P}_j(L)$. $\overline{P}_j$ is normalized to be multiplicative under disjoint union of links and it takes the value 
\[\overline{P}_j(\bigcirc)= a^{-j}q^{j^2} \frac{\qp{a^2q^{2-2j}}{j}}{\qpp{j}}\]
 on the unknot. The reduced $j$-colored HOMFLY-PT invariant $P_j$ of a link $L$ with a $j$-colored component is defined via:
 \[ \overline{P}_j(L)=P_j(L)\overline{P}_j(\bigcirc)\]
 
\begin{remark}\label{rem:HOMFLY-links} The reduced $j$-colored HOMFLY-PT invariants of knots are in fact Laurent polynomials in $a$ and $q$. More generally, the reduced $j$-colored invariant of an $l$-component link is a rational function in $a$ and $q$ with the standard denominator $\qpp{j}^{l-1}$. This explains why we consider different  generating series in Theorems~\ref{thm:knots} and~\ref{thm:links}. For links with $l$ components, one should work with the generating function:
\[\sum_{j\geq 0} \qpp{j}^{l-2} P_j(L) x^j\]
This has the chance of being a specialization of the generating function of a quiver $Q_L$, whose vertices correspond to generators for a suitably defined and finite-dimensional \emph{totally reduced HOMFLY-PT homology} of $L$.
\end{remark}

\begin{figure}[h]
\begin{gather*}
\scalemath{0.75}{
\begin{tikzpicture} [scale=.75,anchorbase]
\draw[thick,directed=.9] (0.5,0) to [out=90,in=270] (-0.5,1.5);
\draw[white, line width=.15cm] (-0.5,0) to [out=90,in=270] (0.5,1.5);
\draw[thick,directed=.9] (-0.5,0) to [out=90,in=270] (0.5,1.5);
\node at (-0.8,0.1) {\scriptsize $k$};
\node at (0.8,0.1) {\scriptsize $l$};
\node at (-0.8,1.4) {\scriptsize $l$};
\node at (0.8,1.4) {\scriptsize $k$};
\end{tikzpicture} \stackrel{k\geq l}{=} 
\sum_{h=0}^{l} (-q)^{h-l}
\begin{tikzpicture} [scale=.75,anchorbase]
\draw[thick, directed=.55] (-0.5,0) to (-0.5,0.5);
\draw[thick, directed=.60] (-0.5,0.5) to (-0.5,1);
\draw[thick, directed=.65] (-0.5,1) to (-0.5,1.5);
\draw[thick, directed=.55] (0.5,0) to (0.5,0.5);
\draw[thick, directed=.60] (0.5,0.5) to (0.5,1);
\draw[thick, directed=.65] (0.5,1) to (0.5,1.5);
\draw[thick, directed=.55] (-0.5,0.9) to (0.5,1.1);
\draw[thick, directed=.55] (0.5,0.4) to (-0.5,.6);
\node at (-0.8,0.1) {\scriptsize $k$};
\node at (0.8,0.1) {\scriptsize $l$};
\node at (-.8,1.4) {\scriptsize $l$};
\node at (.8,1.4) {\scriptsize $k$};
\node at (0,0.25) {\scriptsize $h$};
\node at (0,1.25) {\scriptsize $$};
\end{tikzpicture}
,\quad 
\begin{tikzpicture} [scale=.75,anchorbase]
\draw[thick,directed=.9] (0.5,0) to [out=90,in=270] (-0.5,1.5);
\draw[white, line width=.15cm] (-0.5,0) to [out=90,in=270] (0.5,1.5);
\draw[thick,directed=.9] (-0.5,0) to [out=90,in=270] (0.5,1.5);
\node at (-0.8,0.1) {\scriptsize $k$};
\node at (0.8,0.1) {\scriptsize $l$};
\node at (-0.8,1.4) {\scriptsize $l$};
\node at (0.8,1.4) {\scriptsize $k$};
\end{tikzpicture} \stackrel{k\leq l}{=} 
\sum_{h=0}^{k} (-q)^{h-k}
\begin{tikzpicture} [scale=.75,anchorbase]
\draw[thick, directed=.55] (-0.5,0) to (-0.5,0.5);
\draw[thick, directed=.60] (-0.5,0.5) to (-0.5,1);
\draw[thick, directed=.65] (-0.5,1) to (-0.5,1.5);
\draw[thick, directed=.55] (0.5,0) to (0.5,0.5);
\draw[thick, directed=.60] (0.5,0.5) to (0.5,1);
\draw[thick, directed=.65] (0.5,1) to (0.5,1.5);
\draw[thick, directed=.55] (-0.5,0.4) to (0.5,0.6);
\draw[thick, directed=.55] (0.5,0.9) to (-0.5,1.1);
\node at (-0.8,0.1) {\scriptsize $k$};
\node at (0.8,0.1) {\scriptsize $l$};
\node at (-.8,1.4) {\scriptsize $l$};
\node at (.8,1.4) {\scriptsize $k$};
\node at (0,0.25) {\scriptsize $h$};
\node at (0,1.25) {\scriptsize $$};
\end{tikzpicture}\quad
\boxed{{\text{For negative } \atop \text{crossings: } q \mapsto q^{-1} }}}
\\
\scalemath{0.75}{
\xy
(0,0)*{
\begin{tikzpicture} [scale=1]
\draw[thick, directed=.55] (-0.6,0) to (0,1);
\draw[thick, directed=.55] (0.6,0) to (0.3,0.5);
\draw[thick, directed=.55] (0.3,0.5) to (0,1);
\draw[thick, directed=.55] (0,0) to (0.3,0.5);
\draw[thick, directed=.55] (0,1) to (0,1.5);
\node at (-0.8,0.1) {\scriptsize $k$};
\node at (-0.2,0.1) {\scriptsize $l$};
\node at (0.8,0.1) {\scriptsize $m$};
\node at (0.7,0.8) {\scriptsize $l+m$};
\node at (0.7,1.4) {\scriptsize $k+l+m$};
\end{tikzpicture}
}
\endxy
\;
=
\xy
(0,0)*{
\begin{tikzpicture} [scale=1]
\draw[thick, directed=.55] (0.6,0) to (0,1);
\draw[thick, directed=.55] (-0.6,0) to (-0.3,0.5);
\draw[thick, directed=.55] (-0.3,0.5) to (0,1);
\draw[thick, directed=.55] (0,0) to (-0.3,0.5);
\draw[thick, directed=.55] (0,1) to (0,1.5);
\node at (-0.8,0.1) {\scriptsize $k$};
\node at (0.2,0.1) {\scriptsize $l$};
\node at (0.8,0.1) {\scriptsize $m$};
\node at (-0.7,0.8) {\scriptsize $k+l$};
\node at (-0.7,1.4) {\scriptsize $k+l+m$};
\end{tikzpicture}
}
\endxy,
\quad
\xy
(0,0)*{
\begin{tikzpicture} [scale=1]
\draw[thick, rdirected=.55] (-0.6,0) to (0,-1);
\draw[thick, rdirected=.55] (0.6,0) to (0.3,-0.5);
\draw[thick, rdirected=.55] (0.3,-0.5) to (0,-1);
\draw[thick, rdirected=.55] (0,0) to (0.3,-0.5);
\draw[thick, rdirected=.55] (0,-1) to (0,-1.5);
\node at (-0.8,-0.1) {\scriptsize $k$};
\node at (-0.2,-0.1) {\scriptsize $l$};
\node at (0.8,-0.1) {\scriptsize $m$};
\node at (0.7,-0.8) {\scriptsize $l+m$};
\node at (0.7,-1.4) {\scriptsize $k+l+m$};
\end{tikzpicture}
}
\endxy
\;
=
\xy
(0,0)*{
\begin{tikzpicture} [scale=1]
\draw[thick, rdirected=.55] (0.6,0) to (0,-1);
\draw[thick, rdirected=.55] (-0.6,0) to (-0.3,-0.5);
\draw[thick, rdirected=.55] (-0.3,-0.5) to (0,-1);
\draw[thick, rdirected=.55] (0,0) to (-0.3,-0.5);
\draw[thick, rdirected=.55] (0,-1) to (0,-1.5);
\node at (-0.8,-0.1) {\scriptsize $k$};
\node at (0.2,-0.1) {\scriptsize $l$};
\node at (0.8,-0.1) {\scriptsize $m$};
\node at (-0.7,-0.8) {\scriptsize $k+l$};
\node at (-0.7,-1.4) {\scriptsize $k+l+m$};
\end{tikzpicture}
}
\endxy,
}
\\
\scalemath{0.75}{
\xy
(0,0)*{
\begin{tikzpicture} [scale=1]
\draw[thick, directed=.75] (0,0) to (0,0.4);
\draw[thick, directed=.75] (0,1) to (0,1.4);
\draw[thick, directed=.55] (0,0.4) to [out=0,in=0](0,1);
\draw[thick, directed=.55] (0,0.4) to [out=180,in=180](0,1);
\node at (-0.5,1.3) {\scriptsize $k+l$};
\node at (-0.4,0.7) {\scriptsize $k$};
\node at (0.4,0.7) {\scriptsize $l$};
\node at (-0.5,0.1) {\scriptsize $k+l$};
\end{tikzpicture}
}
\endxy 
=
{k+l \brack l}\;
\xy
(0,0)*{
\begin{tikzpicture} [scale=1]
\draw[thick, directed=.55] (0,0) to (0,1.4);
\node at (0.4,0.1) {\scriptsize $k+l$};
\node[opacity=0] at (0.4,1.3) {\scriptsize $k+l$};
\end{tikzpicture}
}
\endxy 
,\quad
\xy
(0,0)*{
\begin{tikzpicture} [scale=1]
\draw[thick, directed=.55] (-0.5,0) to (-0.5,0.5);
\draw[thick, directed=.60] (-0.5,0.5) to (-0.5,1);
\draw[thick, directed=.65] (-0.5,1) to (-0.5,1.5);
\draw[thick, directed=.55] (0.5,0) to (0.5,0.5);
\draw[thick, directed=.60] (0.5,0.5) to (0.5,1);
\draw[thick, directed=.65] (0.5,1) to (0.5,1.5);
\draw[thick, directed=.55] (-0.5,0.4) to (0.5,0.6);
\draw[thick, directed=.55] (0.5,0.9) to (-0.5,1.1);
\node at (-0.8,0.1) {\scriptsize $k$};
\node at (0.8,0.1) {\scriptsize $l$};
\node at (-1,0.75) {\scriptsize $k-s$};
\node at (1,0.75) {\scriptsize $l+s$};
\node at (-1.2,1.4) {\scriptsize $k-s+r$};
\node at (1.2,1.4) {\scriptsize $l+s-r$};
\node at (0,0.25) {\scriptsize $s$};
\node at (0,1.25) {\scriptsize $r$};
\end{tikzpicture}
}
\endxy
\!\!\!\!\!\!\!= \;
\sum_{t} {k-l+r-s \brack t}
\xy
(0,0)*{
\begin{tikzpicture} [scale=1]
\draw[thick, directed=.55] (-0.5,0) to (-0.5,0.5);
\draw[thick, directed=.60] (-0.5,0.5) to (-0.5,1);
\draw[thick, directed=.65] (-0.5,1) to (-0.5,1.5);
\draw[thick, directed=.55] (0.5,0) to (0.5,0.5);
\draw[thick, directed=.60] (0.5,0.5) to (0.5,1);
\draw[thick, directed=.65] (0.5,1) to (0.5,1.5);
\draw[thick, directed=.55] (0.5,0.4) to (-0.5,0.6);
\draw[thick, directed=.55] (-0.5,0.9) to (0.5,1.1);
\node at (-0.8,0.1) {\scriptsize $k$};
\node at (0.8,0.1) {\scriptsize $l$};
\node at (-1.2,0.75) {\scriptsize $k+r-t$};
\node at (1.2,0.75) {\scriptsize $l-r+t$};
\node at (-1.2,1.4) {\scriptsize $k-s+r$};
\node at (1.2,1.4) {\scriptsize $l+s-r$};
\node at (0,0.25) {\scriptsize $r-t$};
\node at (0,1.25) {\scriptsize $s-t$};
\end{tikzpicture}
}
\endxy
}
\\
\scalemath{0.75}{
\xy
(0,0)*{
\begin{tikzpicture} [scale=1]
\draw[thick, directed=.75] (0,0) to (0,0.4);
\draw[thick, directed=.75] (0,1) to (0,1.4);
\draw[thick, rdirected=.50] (0,0.4) to [out=0,in=0](0,1);
\draw[thick, directed=.55] (0,0.4) to [out=180,in=180](0,1);
\node at (-0.25,1.3) {\scriptsize $k$};
\node at (-0.6,0.7) {\scriptsize $k+l$};
\node at (0.4,0.7) {\scriptsize $l$};
\node at (-0.25,0.1) {\scriptsize $k$};
\end{tikzpicture}
}
\endxy 
=
\xy
(0,0)*{
\begin{tikzpicture} [scale=1]
\draw[thick, directed=.75] (0,0) to (0,0.4);
\draw[thick, directed=.75] (0,1) to (0,1.4);
\draw[thick, directed=.50] (0,0.4) to [out=0,in=0](0,1);
\draw[thick, rdirected=.55] (0,0.4) to [out=180,in=180](0,1);
\node at (-0.25,1.3) {\scriptsize $k$};
\node at (0.6,0.7) {\scriptsize $k+l$};
\node at (-0.4,0.7) {\scriptsize $l$};
\node at (-0.25,0.1) {\scriptsize $k$};
\end{tikzpicture}
}
\endxy 
=
a^{-l}q^{l(l+k)}\frac{\qp{a^2q^{2-2k-2l}}{l}}{\qpp{l}}\;
\xy
(0,0)*{
\begin{tikzpicture} [scale=1]
\draw[thick, directed=.55] (0,0) to (0,1.4);
\node at (0.2,0.1) {\scriptsize $k$};
\node[opacity=0] at (0.2,1.3) {\scriptsize $k$};
\end{tikzpicture}
}
\endxy ,\quad
\xy
(0,0)*{
\begin{tikzpicture} [scale=.75]
\draw[thick, directed=.55] (-0.5,0) to (-0.5,0.5);
\draw[thick, directed=.60] (-0.5,0.5) to (-0.5,1);
\draw[thick, directed=.65] (-0.5,1) to (-0.5,1.5);
\draw[thick, rdirected=.55] (0.5,0) to (0.5,0.5);
\draw[thick, rdirected=.60] (0.5,0.5) to (0.5,1);
\draw[thick, rdirected=.65] (0.5,1) to (0.5,1.5);
\draw[thick, rdirected=.55] (-0.5,0.4) to [out=340,in=200] (0.5,0.4);
\draw[thick, rdirected=.55] (0.5,1.1) to [out=160,in=20] (-0.5,1.1);
\node at (-0.8,0.1) {\scriptsize $1$};
\node at (0.8,0.1) {\scriptsize $1$};
\node at (-.8,0.75) {\scriptsize $2$};
\node at (.8,0.75) {\scriptsize $2$};
\node at (-.8,1.4) {\scriptsize $1$};
\node at (.8,1.4) {\scriptsize $1$};
\node at (0,0.55) {\scriptsize $1$};
\node at (0,.95) {\scriptsize $1$};
\end{tikzpicture}
}
\endxy
= \xy
(0,0)*{
\begin{tikzpicture} [scale=.75]
\draw[thick, directed=.55] (-0.5,0) to (-0.5,0.2) to [out=90,in=90] (.5,.2) to (.5,0);
\draw[thick, directed=.55] (.5,1.5) to (.5,1.3) to [out=270,in=270] (-0.5,1.3) to (-0.5,1.5);
\node at (-0.7,0.1) {\scriptsize $1$};
\node at (0.7,0.1) {\scriptsize $1$};
\node at (-0.7,1.4) {\scriptsize $1$};
\node at (0.7,1.4) {\scriptsize $1$};
\end{tikzpicture}
}
\endxy\!\!\!\!\!
+ 
a^{-1}q^3\frac{\qp{a^2q^{-4}}{1}}{\qpp{1}}
\xy
(0,0)*{
\begin{tikzpicture} [scale=.75]
\draw[thick, directed=.55] (-0.5,0) to (-0.5,1.5);
\draw[thick, rdirected=.55] (0.5,0) to (0.5,1.5);
\node at (-0.7,0.1) {\scriptsize $1$};
\node at (0.7,0.1) {\scriptsize $1$};
\node at (-0.7,1.4) {\scriptsize $1$};
\node at (0.7,1.4) {\scriptsize $1$};
\end{tikzpicture}
}
\endxy
}
\end{gather*}
\caption{The colored HOMFLY-PT skein theory cheat sheet, \cite{CKM,TVW}.}
\label{fig:skein}
\end{figure}
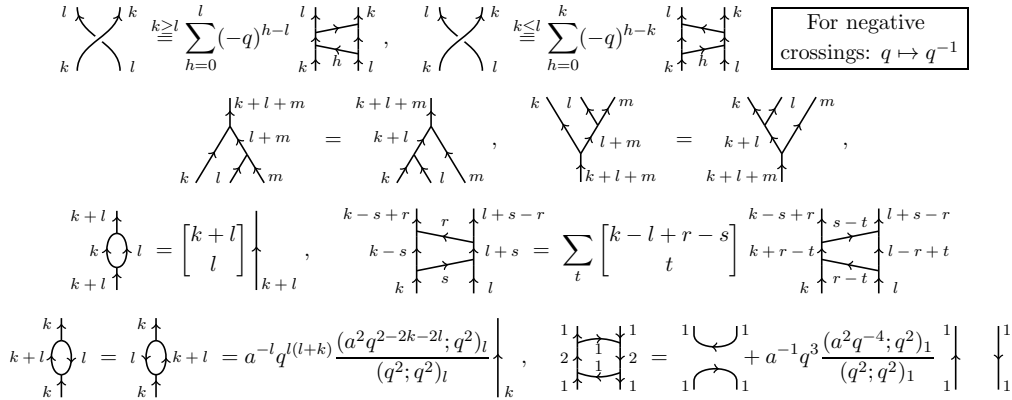

The unreduced colored HOMFLY-PT invariants can be computed via a skein theory of webs that goes back to \cite{MOY}. We summarize the main features of this tool in Figure~\ref{fig:skein}. Given a link or tangle diagram $D$, there are local \textit{crossing rules} that allow the \textit{skein module evaluation} $\langle D \rangle $ to be rewritten as a $\C[a^{\pm 1}](q)$-linear combination of \textit{webs}, which are planar trivalent graphs, each decorated by a non-negative integer flow. Second, there are local linear relations on webs, which allow each closed web to be evaluated to a scalar multiple of the empty web. Thus, each link diagram evaluates in two steps to such a multiple, and the coefficient is the unreduced colored HOMFLY-PT invariant of the link.

\subsubsection{Basic webs}
\label{sec:basicwebs}
The $j$-colored HOMFLY-PT skein invariant $\langle \tau \rangle_j$ of a $4$-ended tangle $\tau$ can be expanded unambiguously as a $\C[a^{\pm 1}](q)$-linear combination of the following basis webs. For each configuration of boundary orientation there are $j+1$ such basis elements indexed by non-negative integers $k\leq j$. 

\[ UP[j,k] \;= \xy
(0,0)*{
\begin{tikzpicture} [scale=.75]
\draw[thick, directed=.55] (-0.5,0) to (-0.5,0.5);
\draw[thick, directed=.60] (-0.5,0.5) to (-0.5,1);
\draw[thick, directed=.65] (-0.5,1) to (-0.5,1.5);
\draw[thick, directed=.55] (0.5,0) to (0.5,0.5);
\draw[thick, directed=.60] (0.5,0.5) to (0.5,1);
\draw[thick, directed=.65] (0.5,1) to (0.5,1.5);
\draw[thick, directed=.55] (0.5,0.4) to (-0.5,0.55);
\draw[thick, directed=.55] (-0.5,0.95) to (0.5,1.1);
\node at (-0.7,0.1) {\scriptsize $j$};
\node at (0.7,0.1) {\scriptsize $j$};
\node at (-0.7,1.4) {\scriptsize $j$};
\node at (0.7,1.4) {\scriptsize $j$};
\node at (0,0.25) {\scriptsize $k$};
\node at (0,1.25) {\scriptsize $k$};
\end{tikzpicture}
}
\endxy 
\;,\quad
OP[j,k] \;= \xy
(0,0)*{
\begin{tikzpicture} [scale=.75]
\draw[thick, directed=.55] (-0.5,0) to (-0.5,0.5);
\draw[thick, directed=.60] (-0.5,0.5) to (-0.5,1);
\draw[thick, directed=.65] (-0.5,1) to (-0.5,1.5);
\draw[thick, rdirected=.55] (0.5,0) to (0.5,0.5);
\draw[thick, rdirected=.60] (0.5,0.5) to (0.5,1);
\draw[thick, rdirected=.65] (0.5,1) to (0.5,1.5);
\draw[thick, directed=.55] (-0.5,0.45) to [out=20,in=160] (0.5,0.45);
\draw[thick, directed=.55] (0.5,1.05) to [out=200,in=340] (-0.5,1.05);
\node at (-0.7,0.1) {\scriptsize $j$};
\node at (0.7,0.1) {\scriptsize $j$};
\node at (-0.7,1.4) {\scriptsize $j$};
\node at (0.7,1.4) {\scriptsize $j$};
\node at (0,1.25) {\scriptsize $k$};
\node at (0,.25) {\scriptsize $k$};
\end{tikzpicture}
}
\endxy
\;,\quad
RI[j,k]\;=
\xy
(0,0)*{
\begin{tikzpicture} [scale=.75]
\draw[thick, directed=.55] (-0.5,0) to (-0.5,0.5);
\draw[thick, rdirected=.60] (-0.5,0.5) to (-0.5,1);
\draw[thick, rdirected=.65] (-0.5,1) to (-0.5,1.5);
\draw[thick, rdirected=.55] (0.5,0) to (0.5,0.5);
\draw[thick, directed=.60] (0.5,0.5) to (0.5,1);
\draw[thick, directed=.65] (0.5,1) to (0.5,1.5);
\draw[thick, directed=.55] (-0.5,0.45) to [out=20,in=160] (0.5,0.45);
\draw[thick, rdirected=.55] (0.5,1.05) to [out=200,in=340] (-0.5,1.05);
\node at (-0.7,0.1) {\scriptsize $j$};
\node at (0.7,0.1) {\scriptsize $j$};
\node at (-0.7,1.4) {\scriptsize $j$};
\node at (0.7,1.4) {\scriptsize $j$};
\node at (0,1.25) {\scriptsize $k$};
\end{tikzpicture}
}
\endxy
\]

\subsubsection{Twist rules}
\label{sec:TR}
Adjoining crossings on the top ($T$) or on the right ($R$) induces linear maps between the skein modules for $4$-ended tangles (with various boundary orientations), which can be written in the above basis as follows.

\begin{enumerate}
\item $TUP[j,k]=\sum_{h=k}^j (-q)^{h-j} q^{k^2} {h \brack k}_+ UP[j,h]$
\item $RUP[j,k]=\sum_{h=0}^k (-q)^{h-j} a^{h-j} q^{-2kh+k^2 +j^2} {j-h \brack k-h}_+ OP[j,h]$
\item $TOP[j,k]=\sum_{h=k}^j (-q)^{h} a^k  q^{k^2-2j k} {h \brack k}_+ RI[j,h] $
\item $ROP[j,k]=\sum_{h=0}^k (-q)^{h-j} a^{k-j} q^{2h(j-k)+  (k-j)^2} {j-h \brack k-h}_+ UP[j,h]$
\item $TRI[j,k]=\sum_{h=k}^j (-q)^{h} a^h  q^{k^2-2jh} {h \brack k}_+ OP[j,h] $
\item $RRI[j,k]=\sum_{h=0}^k (-q)^{h} q^{h(2j-2k) + k^2 -j^2} {j-h \brack k-h}_+ RI[j,h]$
\end{enumerate}
These twist rules were computed in \cite{Wed1} up to an overall monomial factor using skew Howe duality, see also \cite{Wed2}. However, they can also be deduced directly from the local skein relations in Figure~\ref{fig:skein}. Here we have chosen the normalization under which framing changes in a $j$-colored strand act by powers of $f(j)=(-q)^{-j} a^{-j} q^{j^2}$. More specifically, a $j$-colored Reidemeister I tangle with a positive crossing equals the $f(j)$ multiple of the trivial tangle. 

\begin{example}
The twist rule for $TUP[j,0]$ is precisely the positive crossing rule illustrated in Figure~\ref{fig:skein}. 
\end{example}

\subsubsection{Closures}
\label{sec:clos}
It is well known that any closed web evaluates to a scalar in $\C[a^\pm](q)$. Moreover, if the closed web has a $j$-edge, then the scalar is a multiple of the scalar associated to a $j$-unknot. We call the result of dividing the scalar by this unknot the $j$-reduced evaluation. 

There is a unique way of producing closed webs from the basic webs $UP[j,k]$, $RI[j,k]$ by joining end points by planar arcs. We denote the operation which sends such a web to the $j$-reduced evaluation of its closure by $Cl(-)$. For the webs $OP[j,k]$, there are two ways of closing, either connecting boundary points by arcs running North-South or East-West. We denote the two corresponding closure operations by $Cl_{NS}(-)$ and $Cl_{EW}(-)$ respectively.

If $\tau$ is a $j$-colored rational tangle, interpreted as an element of the HOMFLY skein module, then the operation(s) $Cl$ compute the reduced $j$-colored HOMFLY polynomials of the 2-bridge link(s) given by the closure(s) of $\tau$.

The reduced closure satisfies 

\begin{itemize}
\item $Cl ( UP[j,k])=a^{-j} q^{j^2+k^2} 
\frac{\qp{a^2 q^{2-2j-2k}}{j}}{\qpp{j}} {j\brack k}_+ $.
\item $Cl ( RI[j,k])=a^{-j} q^{j^2+(j-k)^2} 
\frac{\qp{a^2 q^{2-4j+2k}}{j}}{\qpp{j}} {j\brack k}_+$ .

\item $Cl_{NS} ( OP[j,k])=a^{-j+k} q^{(j-k)^2} 
\frac{\qp{a^2 q^{2-2j}}{j-k}}{\qpp{j-k}} {j\brack k}_+$.
\item $Cl_{EW} ( OP[j,k])=a^{-k} q^{k^2} 
\frac{\qp{a^2 q^{2-2j}}{k}}{\qpp{k}} {j\brack k}_+$.
\end{itemize}
It is straightforward to deduce these closure rules from the skein relations in Figure~\ref{fig:skein}. For the sake of demonstrating how skein module computations work in practice, we check the first closure rule. 
\begin{align*}
\left\langle\;
\xy
(0,0)*{
\begin{tikzpicture} [scale=.65]
\draw[thick] (-0.5,0) to [out=270,in=270] (-1,0) to (-1,1.5) to [out=90,in=90] (-.5,1.5);
\draw[thick] (0.5,0) to [out=270,in=270] (1,0) to (1,1.5) to [out=90,in=90] (.5,1.5);
\draw[thick, directed=.55] (-0.5,0) to (-0.5,0.5);
\draw[thick, directed=.60] (-0.5,0.5) to (-0.5,1);
\draw[thick] (-0.5,1) to (-0.5,1.5);
\draw[thick, directed=.55] (0.5,0) to (0.5,0.5);
\draw[thick, directed=.60] (0.5,0.5) to (0.5,1);
\draw[thick, directed=.65] (0.5,1) to (0.5,1.5);
\draw[thick, directed=.55] (0.5,0.4) to (-0.5,0.6);
\draw[thick, directed=.55] (-0.5,0.9) to (0.5,1.1);
\node at (-0.7,0.1) {\scriptsize $j$};
\node at (0.7,0.1) {\scriptsize $j$};
\node at (0,0.25) {\scriptsize $k$};
\node at (0,1.25) {\scriptsize $k$};
\end{tikzpicture}
}
\endxy  
\;\right\rangle
&= a^{-j} q^{j(j+k)} 
\frac{\qp{a^2 q^{2-2j-2k}}{j}}{\qpp{j}}
\left\langle
\xy
(0,0)*{
\begin{tikzpicture} [scale=.65]
\draw[thick] (0.5,0) to [out=270,in=270] (1,0) to (1,1.5) to [out=90,in=90] (.5,1.5);
\draw[thick, directed=.60] (-0.5,0.6) to (-0.5,.9);
\draw[thick, directed=.55] (0.5,0) to (0.5,0.5);
\draw[thick, directed=.60] (0.5,0.5) to (0.5,1);
\draw[thick, directed=.65] (0.5,1) to (0.5,1.5);
\draw[thick, directed=.55] (0.5,0.4) to (-0.5,0.6);
\draw[thick, directed=.55] (-0.5,0.9) to (0.5,1.1);
\node at (0.7,0.1) {\scriptsize $j$};
\node at (0,.25) {\scriptsize $k$};
\end{tikzpicture}
}
\endxy  
\;\right\rangle\\
&= a^{-j} q^{j^2+k^2} 
\frac{\qp{a^2 q^{2-2j-2k}}{j}}{\qpp{j}}{j\brack k}_+\!\!\!
\left\langle
\xy
(0,0)*{
\begin{tikzpicture} [scale=.65]
\draw[thick] (.5,1.5) to (0.5,0) to [out=270,in=270] (1,0) to (1,1.5) to [out=90,in=90] (.5,1.5);
\node at (0.7,0.1) {\scriptsize $j$};
\end{tikzpicture}
}
\endxy  
\;\right\rangle
\end{align*}
Here we have used ${j \brack k}=q^{-k(j-k)}{j \brack k}_+$ to translate between the balanced $q$-binomial appearing in Figure~\ref{fig:skein} and the positive $q$-binomial used in the rest of the paper. The reduced closure is obtained by dividing by the invariant of a $j$-colored unknot, i.e. by removing the last factor.

\subsubsection{Formulas for the HOMFLY-PT polynomial}
Using the ingredients described in the last section, we can obtain closed formulas for the reduced $j$-colored HOMFLY-PT polynomials of rational links. Indeed, if the rational number $u/v$ has the continued fraction expansion $[a_1,\dots,a_r]$ with $r$ odd, then we can iteratively compute the skein module element $\langle \tau_{u/v} \rangle_j$ associated to the $j$-colored rational tangle $\tau_{u/v}$ by starting with the invariant of the trivial upwards oriented tangle $\langle \tau_{0/1}\rangle_j=UP[j,0]$ and applying the top and right twist operators to it:
\[\langle \tau_{u/v} \rangle_j  =T^{a_r}R^{a_{r-1}}\cdots R^{a_2}T^{a_1}UP[j,0] \] 
Finally, the reduced $j$-colored HOMFLY-PT polynomial of the rational link $L_{u/v}$ is computed by applying the NS-closure operation:
\[P_j(L_{u/v})=Cl_{NS}(\langle \tau_{u/v} \rangle_j)\]
Using this recipe, we get an expression of the following form for the generating function of the reduced colored HOMFLY-PT polynomials:
\begin{gather*}
\sum_{j\geq 0} P_j(L_{u/v})x^j = \!\!\!\!\!\!\!\!\!\!\!\!\!\!\sum_{\substack{0\leq k^1_1\leq\cdots\leq k^1_{a_1}\\ k^1_{a_1}\geq k^2_{1}\geq\cdots\geq k^2_{a_2}\geq 0\\ \cdots \\k^{r-1}_{a_{r-1}}\leq k^{r}_{1}\leq \cdots\leq k^r_{a_r}\leq j}} \!\!\!\!\!\!\!\!\!\!\!\!\!\!\!(-q)^{\#} a^{\#} q^{\#\#} {k^1_{a_1} \brack \#,\dots,\# }{j-k^2_{1} \brack \#,\dots, \# }\cdots  {k^{r-2}_{a_{r-2}} \brack \#, \dots, \# }\\[-6ex]
{\hspace{6cm}\times {j-k^{r-1}_{1} \brack \#,\dots,\# } {j \brack \#,k^r_{a_r} }\frac{\qp{a^2q^{\#}}{\#}}{\qpp{\#}}x^j}
\end{gather*}
Here we have one summation index per crossing and $\#$ and $\#\#$ stand for linear and quadratic expressions in the summation indices respectively. The length of the continued fraction expansion gives the number of quantum multinomial coefficients appearing in the formula (they arise by telescoping quantum binomial coefficients coming from individual crossings). 
Note however, that the summation indices are not independent, but are subject to ordering constraints determined by the sequence of $T$ and $R$ operators. It is possible to algorithmically rewrite such expressions into quiver form \eqref{eqn:quiverform}, but it is more efficient to perform such rewritings into quiver form at every step of the computation of the skein module element $\langle \tau_{u/v}\rangle_j$. We will describe this in Section~\ref{sec:links}.

\subsection{Example: torus links}

The quiver sum associated to the $T_{2,n}$ torus link is computed as follows.
\begin{align*}
Cl (T^n  UP[j,0])=
&\sum_{0\le k_1\le k_2\le\cdots\le k_n\le j}  \!\!\!\!\!\!\!(-q)^{\sum_{i=1}^{n} k_i -n j} a^{-j} q^{j^2 +\sum_{i=1}^{n} k_i^2}{k_2 \brack k_1}_+\cdots{k_n \brack k_{n-1}}_+   \\ &\hspace{5cm}\times {j \brack k_n}_+ \frac{\qp{a^2q^{2-2j-2k_n}}{j}}{\qpp{j}}
\\
&\!\!\!\!\!\!\!\!\!\!\!\!\!\!\!\!\!\!\!\!\!\!\!\!\!\!\!\!=\sum_{0\le k_1\le k_2\le\cdots\le k_n\le j}  \frac{(-q)^{\sum_{i=1}^{n} k_i -n j} a^{-j} q^{j^2 +\sum_{i=1}^{n} k_i^2}  \qp{a^2q^{2-2j-2k_n}}{j}}{\qpp{k_1}\qpp{k_2-k_1}\cdots \qpp{k_n-k_{n-1}}\qpp{j-k_n}}
\end{align*}
After introducing the new variables $d_1=k_1$, $d_i=k_i-k_{i-1}$ for $2\leq i \leq n$ and $d_{n+1}=j-k_n$, this is sum is almost in quiver form. One application of Lemma~\ref{lemlong} allows to rewrite it in quiver form with $2n+2$ summation indices. Later we will see that for torus knots of this form, i.e. when $n$ is odd, the quiver sum can be rewritten as a polynomial quiver sum with $n$ summation indices.

\section{Quivers and rational tangles}
\label{sec:links}
In this section we define what it means for the skein module element of a rational tangle to be encoded in quiver form. We then argue that the trivial tangle has such an expression and that the operations of adding top- and right twists preserve this property. This gives an iterative way of computing a quiver form for the generating function of the colored HOMFLY-PT polynomials of rational links. 

\begin{definition}
Let $\tau$ be a $4$-ended tangle with boundary orientations of type $X$ and let 
$\langle \tau \rangle_j = \sum_{k=0}^j C_k X[j,k]$ be
its $j$-colored skein module evaluation, expanded in terms of the basis webs from Section~\ref{sec:basicwebs}. Then we define the \textit{rescaled skein module evaluation} of the $j$-colored tangle $\tau$ to be 
\[\langle \tau \rangle_j^\prime = \sum_{k=0}^j C_k {j \brack k}_+ X[j,k]\]
and collect all of these into a generating function $P^\prime(\tau)=\sum_{j\geq 0} \langle \tau \rangle_j^\prime$.
\end{definition}
Note that this rescaling is invertible and, thus, $\langle \tau \rangle_j$ is determined by $\langle \tau \rangle_j^\prime$.

\begin{definition}\label{def:equiverform} A generating function in \textit{quiver form} for skein module elements of type $X$ is an expression of the form
\begin{equation}
\label{eqn:equiverform}
\sum_{\textbf{d}=(d_1,\dots,d_{m+n})\in \N^{m+n}} \!\!\!\!(-q)^{S\cdot\textbf{d}} a^{A\cdot\textbf{d}} q^{\textbf{d}\cdot Q\cdot \textbf{d}^t} {d_1+\cdots+d_{m+n} \brack d_1,\ldots,d_{m+n}}  X[d_1+\cdots+d_{m+n},d_1+\cdots+d_m].
\end{equation}
We call the summation indices $d_1,\dots,d_m$ \textit{active} and the remaining $d_{m+1},\dots,d_{m+n}$ \textit{inactive}. As before, $S$ and $A$ are integer vectors of length $m+n$ and $Q$ is an symmetric integer matrix of size $m+n$.

We say a generating function is in \textit{almost quiver form}, if it has an additional $q$-Pochhammer symbol $\qpp{d_{i_1}+\cdots+d_{i_k}}$ as factor in each summand:
\begin{gather}
\label{eqn:equiverform2}
\sum_{\textbf{d}=(d_1,\dots,d_{m+n})\in \N^{m+n}} (-q)^{S\cdot\textbf{d}} a^{A\cdot\textbf{d}} q^{\textbf{d}\cdot Q\cdot \textbf{d}^t} \qpp{K \cdot \textbf{d}} {d_1+\cdots+d_{m+n} \brack d_1,\ldots,d_{m+n}}\\
\hspace{5cm}\nonumber\times  X[d_1+\cdots+d_{m+n},d_1+\cdots+d_m].
\end{gather}
 Here we write $K$ for the vector of length $m+n$ with entries $1$ at the indices $i_j$ for $1 \leq j \leq k$ and zeros elsewhere, so that $d_{i_1}+\cdots+d_{i_k}=K\cdot \textbf{d}$.   
\end{definition}

\begin{notation}\label{not:triple} Later, we will need a more compact way of writing expressions of the form \eqref{eqn:equiverform2} (and its special case \eqref{eqn:equiverform}). To this end, we split $K$, $S$ and $A$ into vectors $K_+$ and $K_-$, $S_+$ and $S_-$ as well as $A_+$ and $A_-$ according to whether they contribute coefficients to active or inactive summation indices. The matrix $Q$ can be written as a $2\times 2$ block matrix with blocks $Q_{++}$, $Q_{+-}$, $Q_{-+}$ and $Q_{--}$ accordingly. This data of a generating function in almost quiver form for skein module elements of type $X$ is then encoded as the following triple:
\begin{equation}
\label{eqn:triple}
  \scalemath{0.75}{ X, \left( \begin{array}{c|c|c} K_+ & S_+ & A_+ \\ \hline\hline K_- & S_- & A_- 
\end{array} \right)},
 \scalemath{0.75}{
 \left( 
\begin{array}{c|c} Q_{++} & Q_{+-} \\ \hline Q_{-+} & Q_{--} 
\end{array} \right)} 
\end{equation} 
The first entry gives the object type. The first matrix encodes the additional $q$-Pochhammer symbol in the first column and the linear exponent of $-q$ and $a$ in the second and third column respectively. The horizontal double line separates coefficients of active summation indices (above) from those of inactive ones (below). The second matrix describes the quadratic part of the exponent of $q$.  
\end{notation}

\begin{theorem} \label{thm:TRonQS}
 The generating functions for the rescaled skein module evaluations $\sum_{j}\langle \tau_{u/v}\rangle^\prime_j$ of a rational tangle $\tau_{u/v}$ can be written in quiver form \eqref{eqn:equiverform}.
\end{theorem}
\begin{proof} The proof proceeds by induction on the number of crossings of the rational tangle. For the trivial tangle the claim is true:
\[\sum_{j\geq 0} {j \brack 0}_{+} UP[j,0]= \sum_{j\geq 0} {j \brack 0,j} UP[j,0]\]
Now suppose that we have brought the rescaled skein module evaluation of $\tau$ into quiver form. Then the unrescaled evaluation has the following expression:
\begin{gather}
\label{eqn:skeinel}
\sum_{j\geq 0}\langle \tau\rangle_j =\!\!\!\!\!\!\!\!\! \sum_{d_1,\dots,d_{m+n}\geq 0}\!\!\!\!\!\!\!\!\! (-q)^{\#} a^{\#} q^{\#\#} {d_{\underline{m}} \brack d_1,\ldots,d_m} {d_{\underline{m+n}}-d_{\underline{m}} \brack d_{m+1},\ldots,d_{m+n}} X[d_{\underline{m+n}},d_{\underline{m}}] 
\end{gather}
where $\#$ denote some linear expressions in the $d_i$ and $\#\#$ denotes some purely quadratic expression in the $d_i$. Additionally, we have abbreviated $d_{\underline{m}}:=d_1+\cdots+d_m$ and $d_{\underline{m+n}}:=d_1+\cdots+d_{m+n}=j$.

Let $T\tau$ denote the tangle obtained by adding a top twist to $\tau$. Its skein module evaluation can be computed by applying the appropriate top twist rule and we get an expression of the following form for $\sum_{j\geq 0}\langle T\tau \rangle_j$:
\begin{gather}
\label{eqn:aftertop} 
\sum_{\substack{d_1,\dots,d_{m+n}\geq 0\\ d_{\underline{m}}\leq h \leq d_{\underline{m+n}} }} (-q)^{\#} a^{\#} q^{\#\#} {d_{\underline{m}} \brack d_1,\ldots,d_m} {d_{\underline{m+n}}-d_{\underline{m}} \brack d_{m+1},\ldots,d_{m+n}} {h \brack d_{\underline{m}} }_+  TX[d_{\underline{m+n}},h] 
\end{gather}
Here $TX$ denotes the type of object obtained from $X$ by applying the top twist.

Now we need to check that the rescaled skein module evaluation can be rewritten in quiver form. To this end, we use Lemma~\ref{lemvlong} to rewrite the second $q$-multinomial in \eqref{eqn:aftertop} as 

\begin{gather*}
{(d_{\underline{m+n}}-h)+(h-d_{\underline{m})} \brack d_{m+1},\ldots,d_{m+n}} 
= \sum_{\substack{b_1+c_1=d_{m+1}\\ \cdots \\ b_n+c_n=d_{m+n} }} q^{\#\#} {d_{\underline{m+n}}-h \brack b_{1},\ldots,b_{n}} {h-d_{\underline{m}} \brack c_{1},\ldots,c_{n}}.
\end{gather*}

Using this and rewriting $h=d_{\underline{m}} + c_{\underline{n}}$ and $d_{\underline{m+n}}= d_{\underline{m}}+ b_{\underline{n}}+ c_{\underline{n}} $, we encounter the following $q$-multinomials in  \eqref{eqn:aftertop}:

\begin{align*}
 {d_{\underline{m}} \brack d_1,\ldots,d_m}
 {c_{\underline{n}} \brack c_{1},\ldots,c_{n}}
 {d_{\underline{m}}+c_{\underline{n}} \brack d_{\underline{m}},c_{\underline{n}}} {b_{\underline{n}} \brack b_{1},\ldots,b_{n}}  \\=
{d_{\underline{m}}+c_{\underline{n}} \brack d_1,\ldots,d_m,c_1,\ldots,c_n} {b_{\underline{n}} \brack b_{1},\ldots,b_{n}}
\end{align*}

After the rescaling multiplication by ${d_{\underline{m+n}} \brack h}_+= {d_{\underline{m}}+c_{\underline{n}}+b_{\underline{n}} \brack d_{\underline{m}}+c_{\underline{n}}, b_{\underline{n}}}$, we obtain the desired expression in quiver form for the generating function $\sum_{j\geq 0} \langle T\tau\rangle_j^\prime$:

\begin{gather*}
\sum_{\substack{d_1,\dots,d_{m}\geq 0\\c_1,\dots,c_{n}\geq 0\\b_1,\dots,b_{n}\geq 0}} \!\!\!\!\!\!(-q)^{\#} a^{\#} q^{\#\#} {d_{\underline{m}}+c_{\underline{n}}+b_{\underline{n}} \brack d_1,\ldots,d_m,c_1,\ldots,c_n,b_1,\ldots,b_n}   TX[d_{\underline{m}}\!+\!c_{\underline{n}}\!+\!b_{\underline{n}},d_{\underline{m}}\!+\!c_{\underline{n}}] 
\end{gather*}
Note for later use, that this top twist has increased the number of active indices by the number of previously inactive indices, while the number of inactive indices remained constant.

The case of right twists is analogous, with the main difference that each active summation index (rather than each inactive one) is split into two indices, one of which remains active, while the other one becomes inactive. We only sketch the analogous computation as far as the $q$-multinomial coefficients are concerned. In the case of the right twist, the new summation variable $h$ is constrained to be between $0$ and $d_{\underline{m}}$. Thus, Lemma~\ref{lemvlong} allows us to split the other $q$-multinomial coefficient in \eqref{eqn:skeinel} as follows:
\[{d_{\underline{m}} \brack d_1,\ldots,d_m}= 
{h + (d_{\underline{m}}-h) \brack d_1,\ldots,d_m} = 
\sum_{\substack{b_1+c_1=d_{1}\\ \cdots \\ b_m+c_m=d_{m} }} q^{\#\#} {h \brack b_{1},\ldots,b_{m}} {d_{\underline{m}}-h \brack c_{1},\ldots,c_{m}}
\]
Since the result of the right twist is a linear combination of objects $RX[j,h]$, the summation indices $b_1,\ldots,b_m$ are the new active indices. The right twist furthermore introduces an additional $q$-binomial and finally, to compare with the claim, we need to rescale by another $q$-binomial. Rewriting the product of these $q$-multinomials in the new variables, using $d_{\underline{n}}:=d_{m+1}+\cdots + d_{m+n}$, we get:

\begin{gather*}
 {b_{\underline{m}} \brack b_1,\ldots,b_m}
 {c_{\underline{m}} \brack c_{1},\ldots,c_{m}}
 {d_{\underline{n}} \brack d_{m+1},\ldots,d_{m+n}} 
 {c_{\underline{m}}+d_{\underline{n}} \brack c_{\underline{m}}, d_{\underline{n}}}   
 {b_{\underline{m}}+c_{\underline{m}}+d_{\underline{n}} \brack b_{\underline{m}},c_{\underline{m}}+ d_{\underline{n}}}   
  \\=
{b_{\underline{m}}+c_{\underline{m}}+d_{\underline{n}} \brack b_1,\ldots,b_m,c_1,\ldots,c_m,d_{m+1},\ldots,d_{m+n}}
\end{gather*}

This proves the claim for right twists and completes the induction step.
\end{proof}

For later reference, we record the action of the twist rules on such generating functions in full detail.  
\begin{proposition}
\label{prop:twisttriple}
Suppose that the generating function for the rescaled skein module evaluations of a tangle has been written in (almost) quiver form and encoded by the triple \eqref{eqn:triple}. Then the results of top twists $T$ and right twists $R$ depend on $X$ as follows:
\begin{gather*}
\scalemath{0.75}{ TUP, \left( \begin{array}{c|c|c} K_+ & S_+ & A_+ \\ \hline K_- & S_- & A_-\\ \hline\hline K_- & S_- -1 & A_- 
\end{array} \right)},
 \scalemath{0.75}{
 \left( 
\begin{array}{c|c|c} Q_{++}+1 & Q_{+-} & Q_{+-} \\ \hline Q_{-+} & Q_{--}& Q_{--}+L \\ \hline Q_{-+} & Q_{--}+U& Q_{--} 
\end{array} \right)},\!\!
\\
\scalemath{0.75}{ RUP, \left( \begin{array}{c|c|c} K_+ & S_+ & A_+ \\ \hline\hline K_+ & S_+ -1 & A_+ -1 \\ \hline K_- & S_- -1& A_- -1 
\end{array} \right)},
 \scalemath{0.75}{
 \left( 
\begin{array}{c|c|c} Q_{++} & Q_{++}+1+L & Q_{+-}+1 \\ \hline Q_{++}+1+U & Q_{++}+2 & Q_{+-}+1 \\ \hline Q_{-+}+1 & Q_{-+}+1 & Q_{--}+1 
\end{array} \right)}
\\
 \scalemath{0.75}{ TOP, \left( \begin{array}{c|c|c} K_+ & S_+ +1 & A_+ +1 \\ \hline K_- & S_- +1 & A_-\\ \hline\hline K_- & S_- & A_-  
\end{array} \right)},
 \scalemath{0.75}{
 \left( 
\begin{array}{c|c|c} Q_{++}-1 & Q_{+-}-1 & Q_{+-}-1 \\ \hline Q_{-+}-1 & Q_{--}& Q_{--}+L \\ \hline Q_{-+}-1 & Q_{--}+U& Q_{--} 
\end{array} \right)},\!\! 
\\
\scalemath{0.75}{ ROP, \left( \begin{array}{c|c|c} K_+ & S_+ & A_+ \\ \hline\hline K_+ & S_+ -1 & A_+ \\ \hline K_- & S_- -1 & A_- -1 
\end{array} \right)},
 \scalemath{0.75}{
 \left( 
\begin{array}{c|c|c} Q_{++} & Q_{++}+L & Q_{+-}+1 \\ \hline Q_{++}+U & Q_{++} & Q_{+-} \\ \hline Q_{-+}+1 & Q_{-+} & Q_{--}+1 
\end{array} \right)} 
\\
 \scalemath{0.75}{ TRI, \left( \begin{array}{c|c|c} K_+ & S_+ +1 & A_+ +1 \\ \hline K_- & S_- +1 & A_- +1 \\ \hline\hline K_- & S_- & A_- 
\end{array} \right)},
 \scalemath{0.75}{
 \left( 
\begin{array}{c|c|c} Q_{++} -1 & Q_{+-} -2 & Q_{+-} -1 \\ \hline Q_{-+} -2 & Q_{--}-2& Q_{--}-1+L \\ \hline Q_{-+}-1 & Q_{--}-1+U& Q_{--} 
\end{array} \right)},\!\!
\\
\scalemath{0.75}{ RRI, \left( \begin{array}{c|c|c} K_+ & S_+ +1 & A_+ \\ \hline\hline K_+ & S_+ & A_+ \\ \hline K_- & S_- & A_- 
\end{array} \right)},
 \scalemath{0.75}{
 \left( 
\begin{array}{c|c|c} Q_{++} & Q_{++}+L & Q_{+-} \\ \hline Q_{++}+U & Q_{++} & Q_{+-}-1 \\ \hline Q_{-+} & Q_{-+}-1 & Q_{--}-1 
\end{array} \right)}
\end{gather*}

Here $1$ denotes matrices with all entries equal to one, whereas $L$ (resp. $U$) denotes square matrices with all entries below (resp. above) the diagonal equal to one and otherwise zero.
\end{proposition}
\begin{proof}
We will only consider the case of the top twist $T$ applied to a tangle of type $UP$, since the other five cases are completely analogous. Recall that the top twist rule is 
$TUP[j,k]=\sum_{h=k}^j (-q)^{h-j} q^{k^2} {h \brack k}_+ UP[j,h]$. In the proof of Theorem~\ref{thm:TRonQS} we have seen that the top twist operation splits each of the inactive summation indices into a new active and a new inactive index, and that the rescaling enables the $q$-multinomials to transform as desired. A na\"{i}ve splitting of the summation indices would correspond to a simple doubling of the corresponding matrix rows (and columns) of our starting triple:
\[
  \scalemath{0.75}{ UP, \left( \begin{array}{c|c|c} K_+ & S_+ & A_+ \\ \hline\hline K_- & S_- & A_- 
\end{array} \right)},
 \scalemath{0.75}{
 \left( 
\begin{array}{c|c} Q_{++} & Q_{+-} \\ \hline Q_{-+} & Q_{--} 
\end{array} \right)} \;\;
\xrightarrow{} \;\;
\scalemath{0.75}{UP, \left( \begin{array}{c|c|c} K_+ & S_+ & A_+ \\ \hline K_- & S_- & A_-\\ \hline\hline K_- & S_-  & A_- 
\end{array} \right)},
 \scalemath{0.75}{
 \left( 
\begin{array}{c|c|c} Q_{++} & Q_{+-} & Q_{+-} \\ \hline Q_{-+} & Q_{--}& Q_{--} \\ \hline Q_{-+} & Q_{--}& Q_{--} 
\end{array} \right)}
\]
However, the splitting accomplished by Lemma~\ref{lemvlong} additionally introduces the correction factor \eqref{eqn:triangcorr}, which appears as the upper and lower triangular block matrices $U$ and $L$ in:
\[ \scalemath{0.75}{ UP, \left( \begin{array}{c|c|c} K_+ & S_+ & A_+ \\ \hline K_- & S_- & A_-\\ \hline\hline K_- & S_- -1 & A_- 
\end{array} \right)},
 \scalemath{0.75}{
 \left( 
\begin{array}{c|c|c} Q_{++}+1 & Q_{+-} & Q_{+-} \\ \hline Q_{-+} & Q_{--}& Q_{--}+L \\ \hline Q_{-+} & Q_{--}+U& Q_{--} 
\end{array} \right)}
\] The remaining corrections relative to the na\"{i}ve doubling are all due to monomial factors in the twist rule. The factor $(-q)^{h-j}$ has the effect of decreasing each entry of the inactive part of the new $S$ vector by one, and the factor $q^{k^2}$ increases each entry of the first block of the old $Q$ matrix by one.
\end{proof}

Now we can assemble an algorithm to find quivers for rational links.

\begin{theorem}\label{thm:links-full} The generating functions in quiver form for the colored HOMFLY-PT polynomials of rational links $L_{u/v}$ can be computed by the following algorithm. 

\begin{enumerate}
\item Start with the (rescaled) skein module evaluation of the trivial tangle
\[P^\prime(\tau_{0/1}) = \sum_{j\geq 0} {j \brack 0} UP[j,0]\]
which is already in quiver form.

\item Act with a twist operation $T$ or $R$ on the generating function of rescaled skein module evaluations $P^\prime(\tau)$ in quiver form as described in the proof of Theorem~\ref{thm:TRonQS}.
\item Repeat step 2 for every crossing as determined by the odd-length continued fraction expansion of the rational number $u/v$. This computes an expression in quiver form \eqref{eqn:equiverform} for $P^\prime(\tau_{u/v})$.
\item Apply the closure rules to close off the tangle. Since we are working with generating functions for rescaled skein module evaluations, the rules from Section~\ref{sec:clos} take the following form:
\begin{itemize}
\item $UP[j,k]\mapsto a^{-j} q^{j^2+k^2} 
\frac{\qp{a^2 q^{2-2j-2k}}{j}}{\qpp{j}}$
\item $OP[j,k] \mapsto a^{-j+k} q^{(j-k)^2} 
\frac{\qp{a^2 q^{2-2j}}{j-k}}{\qpp{j-k}}$
\end{itemize}
Here we write $k$ for the sum of all active summation indices and $j$ for the sum of all summation indices. Now recall that the summands of the resulting generating functions contain a $q$-multinomial with numerator $\qpp{j}$, which we use to cancel the denominators $\qpp{j}$ or $\qpp{j-k}$ that appear in the closure rules. Once this is done, the generating function is in almost quiver form. That means, it would be in quiver form, were it not for the additional factors $\qp{a^2 q^{2-2j-2k}}{j}$ or $\qp{a^2 q^{2-2j}}{j-k}\qp{q^{2+2j-2k}}{k}$. Finally, we use Lemma~\ref{lemlong} to bring the result into quiver form at the expense of doubling the number of summation indices.
\end{enumerate}
\end{theorem}

The following theorem gives a slight variation of the previous algorithm, which we will need later.

\begin{theorem}(\textbf{Product form})\label{thm:prodform} The following two steps are alternatives to the corresponding steps in the algorithm in Theorem~\ref{thm:links-full}.

\begin{enumerate}
\item[2.] Act with a twist operation $T^\prime$ or $R^\prime$ on the generating function $P^\prime(\tau)$ in quiver form by multiplying by a certain $q$-Pochhammer symbol, depending on the twist and object type:
\begin{itemize}
\item $T^\prime UP[j,k]: (-q)^{k-j}q^{k^2}\qp{q^{2+2k}}{j-k}$
\item $T^\prime OP[j,k]: (-q)^{k}a^k q^{k^2-2jk}\qp{q^{2+2k}}{j-k}$
\item $T^\prime RI[j,k]: (-q)^{k}a^k q^{k^2-2jk}\qp{aq^{2+2k-2j}}{j-k}$
\item $R^\prime UP[j,k]: (-q)^{-j}a^{-j}q^{j^2}\qp{a q^{2-2k}}{k}$
\item $R^\prime OP[j,k]: (-q)^{-j}a^{k-j}q^{j^2-2jk}\qp{q^{2+2j-2k}}{k}$
\item $R^\prime RI[j,k]: q^{-j^2}\qp{q^{2+2j-2k}}{k}$
\end{itemize} 

Then use Lemma~\ref{lemlong} to rewrite the result again in quiver form, thereby splitting some summation indices. Finally, replace the objects $X[j,k]$ by $T^\prime X[j,h]$ or $R^\prime X[j,h]$ as appropriate, where $h$ denotes the sum of the new active summation indices.

\item[4.] Finally, close off the tangle by replacing:
\begin{itemize}
\item $UP[j,k]\mapsto a^{-j} q^{j^2} 
\frac{\qp{a^2 q^{2-2j-2k}}{j}}{\qpp{j}}$
\item $OP[j,k] \mapsto a^{-j+k} q^{j^2-2jk} 
\frac{\qp{a^2 q^{2-2j}}{j-k}}{\qpp{j-k}}$
\end{itemize}
and using cancellation and Lemma~\ref{lemlong} to bring the result into quiver form.
\end{enumerate}
\end{theorem}
\begin{proof}
The proof is by straightforward comparison of the action of the twist rules, as explained in the proof of Theorem~\ref{thm:TRonQS}, with the action of multiplying by the listed Pochhammer symbols and then simplifying using  Lemma~\ref{lemlong}. The only subtlety is that we have rescaled each twist rule $T^\prime$ and $R^\prime$ by $q^{h^2-k^2}$ relative to the versions $T$ and $R$ shown in Section~\ref{sec:TR}. Correspondingly, the closure rules are rescaled by $q^{-k^2}$ relative to the versions in Section~\ref{sec:clos}.

We only illustrate the details in the case of $T^\prime UP$, since the others are completely analogous. First note that multiplication by $(-q)^{k-j}q^{k^2}$ transforms our input generating function, encode in triple notation \eqref{eqn:triple}, as follows:
\[  \scalemath{0.75}{ UP, \left( \begin{array}{c|c|c} 0 & S_+ & A_+ \\ \hline\hline 0 & S_- & A_- 
\end{array} \right)},
 \scalemath{0.75}{
 \left( 
\begin{array}{c|c} Q_{++}  & Q_{+-} \\ \hline Q_{-+} & Q_{--} 
\end{array} \right)} 
\xrightarrow{(-q)^{k-j}q^{k^2}} \scalemath{0.75}{ UP, \left( \begin{array}{c|c|c} 0 & S_+ & A_+ \\ \hline\hline 0 & S_- -1 & A_- 
\end{array} \right)},
 \scalemath{0.75}{
 \left( 
\begin{array}{c|c} Q_{++} + 1 & Q_{+-} \\ \hline Q_{-+} & Q_{--} 
\end{array} \right)}
\]
Next we multiply by the $q$-Pochhammer symbol $\qp{q^{2+2k}}{j-k}$ and rewrite the generating function using Lemma~\ref{lemlong}. During this process, we have to split every one of the $n$ inactive summation index $d_i$ into an active index $\alpha_i$ and and inactive index $\beta_i$, and multiply each summand by a monomial $(-q^{2+2k}q^{-1})^{\sum \alpha_i} q^{\sum \alpha_i^2 + 2\sum_{i=1}^{n-1}\alpha_{i+1}(d_1+\cdots+d_i)}$, which we rewrite as $(-q)^{\sum \alpha_i} q^{2k\sum \alpha_i} q^{(\sum \alpha_i)^2} q^{ 2\sum_{i=1}^{n-1}\alpha_{i+1}(\beta_1+\cdots+\beta_i)}$ using $d_i=\alpha_i+\beta_i$. In triple notation we get:
\[\scalemath{0.75}{ UP, \left( \begin{array}{c|c|c} K_+ & S_+ & A_+ \\ \hline K_- & S_- & A_-\\ \hline\hline K_- & S_- -1 & A_- 
\end{array} \right)},
 \scalemath{0.75}{
 \left( 
\begin{array}{c|c|c} Q_{++}+1 & Q_{+-}+1 & Q_{+-} \\ \hline Q_{-+}+1 & Q_{--}+1& Q_{--}+L \\ \hline Q_{-+} & Q_{--}+U& Q_{--} 
\end{array} \right)}\]
Here, the first block contains the old active summation indices summing to $k$, the second block contains the new active indices $\alpha_i$ and the third block contains the $\beta_i$. The $S$ vector and the $Q$ matrix differ from the na\"{i}ve doubling of matrix rows and columns by shifts induced by the correction monomial.

As claimed, this result agrees with the $TUP$ twist rule one from Proposition~\ref{prop:twisttriple} up to a rescaling by $q^{h^2-k^2}$.
\end{proof}

\begin{remark} The avid reader may have noticed that these algorithms make use of the Lemmas~\ref{lemlong} and~\ref{lemvlong}, which implicitly depend on a choice of ordering of the summation indices. As a consequence, the quivers $Q_L$ (more precisely: the numbers of arrows between distinct vertices) depend in a mild way on an ordering of the vertices. We choose a distinguished ordering which arises naturally from the building process of the rational tangle and which allows certain non-trivial re-summations that are essential for the proof of Theorem~\ref{thm:knots}. More details about this appear in Section~\ref{sec:prooflemmas}.   
\end{remark}

\section{Quivers and rational linksge}
\label{sec:knots}
In the introduction we have stated Theorems~\ref{thm:knots} and~\ref{thm:links} which specify the sizes of the quivers that should be associated to knots and links. These are generally smaller than the quivers that are obtained by Theorem~\ref{thm:links-full}. In this section we will prove the following reformulations:
\begin{theorem}\label{thm:knotty}
The generating functions for the colored HOMFLY-PT polynomials of rational knots $K_{u/v}$ can be brought into \textit{polynomial quiver form}, see Definition~\ref{def:knotty}, with $u$ summation indices. This means, there are $S,A\in \Z^u$ and a $u\times u$ symmetric integer matrix $Q_{K_{u/v}}$, such that 
\begin{equation}
\label{eqn:knotty}
\sum_{j\geq 0} P_j(K_{u/v}) x^j =\sum_{\textbf{d}=(d_1,\dots,d_u)\in \N^u} (-q)^{S\cdot\textbf{d}} a^{A\cdot\textbf{d}} q^{\textbf{d}\cdot Q_{K_{u/v}}\cdot\textbf{d}^t} x^{d_1+\cdots+d_u} {d_1+\cdots+d_u \brack d_1,\ldots,d_u}.
\end{equation}

\end{theorem}
\begin{theorem}\label{thm:linky}
The generating functions for the colored HOMFLY-PT polynomials of rational links $L_{u/v}$ can be brought into \textit{quiver form}, see Definition~\ref{def:linky}, with $2u$ summation indices. This means, there are $S,A \in \Z^{2u}$ and a $2u\times 2u$ symmetric integer matrix $Q_{L_{u/v}}$, such that 
\begin{equation}\sum_{j\geq 0} P_j(L_{u/v}) x^j= \sum_{\textbf{d}=(d_1,\dots,d_{2u})\in \N^{2u}} \frac{(-q)^{S\cdot\textbf{d}} a^{A\cdot\textbf{d}} q^{\textbf{d}\cdot Q_{L_{u/v}}\cdot\textbf{d}^t}} {\qpp{d_1}\cdots \qpp{d_{2u}}}{x^{d_1+\cdots+d_{2u}}}.
\end{equation}
\end{theorem}

In order to prove Theorem~\ref{thm:knotty}, we could perform re-summations in the generating functions provided by Theorem~\ref{thm:links-full} to turn an excess number of summation indices into an additional $q$-Pochhammer factor $\qpp{j}$ in each summand. However, it is easier to describe and perform these re-summations in parallel with the construction process of the rational tangle, proceeding along the continued fraction expansion.

\begin{definition}
The set $\cal{F}_{a}(\alpha,\beta,X)$ is defined to be the set of generating functions of skein module elements of type $X$, which can be written in almost quiver form \eqref{eqn:equiverform2} with $\alpha$ active summation indices, $\beta$ inactive summation indices, and an extra $q$-Pochhammer symbol $\qpp{k}$ of length given by the sum of active summation indices. Analogously, $\cal{F}_{i}(\alpha,\beta,X)$ is defined as the set of such generating functions with extra $q$-Pochhammer symbol $\qpp{j-k}$ of length given by the sum of the inactive summation indices. In triple notation, such generating functions are expressed as:
\[
  \scalemath{0.75}{ X, \left( \begin{array}{c|c|c} 1 & S_+ & A_+ \\ \hline\hline 0 & S_- & A_- 
\end{array} \right)},
 \scalemath{0.75}{
 \left( 
\begin{array}{c|c} Q_{++} & Q_{+-} \\ \hline Q_{-+} & Q_{--} 
\end{array} \right)} 
\;\in\;
\cal{F}_{a}(\alpha,\beta,X)
\quad,\quad
\scalemath{0.75}{ X, \left( \begin{array}{c|c|c} 0 & S_+ & A_+ \\ \hline\hline 1 & S_- & A_- 
\end{array} \right)},
 \scalemath{0.75}{
 \left( 
\begin{array}{c|c} Q_{++} & Q_{+-} \\ \hline Q_{-+} & Q_{--} 
\end{array} \right)} 
\;\in \;
\cal{F}_{i}(\alpha,\beta,X)
\]
\end{definition}

Our first goal is to prove the following proposition about generating functions of rescaled skein module evaluations of rational tangles that close to knots after adding a top twist.

\begin{proposition}\label{prop:newsumalgo} If $u$ is odd, then $P^\prime(\tau_{(u-v)/v}) \in \cal{F}_{a}((u-v)/2,v,UP)$ or $\cal{F}_{i}(u-v,v/2,RI)$ depending on whether $v$ is odd or even.
\end{proposition}
The proof uses the following three lemmas.

\begin{lemma} \label{lem:poch1} Squares of top and right twists act as follows on generating functions in almost quiver form:
\begin{itemize}
\item[(A)] If $P^\prime(\tau) \in \cal{F}_{a}(\alpha,\beta,UP)$, then $P^\prime(T^2\tau) \in \cal{F}_{a}(\alpha+\beta,\beta,UP)$
\item[(A)] If $P^\prime(\tau) \in \cal{F}_{a}(\alpha,\beta,OP)$, then $P^\prime(T^2\tau) \in \cal{F}_{a}(\alpha+\beta,\beta,OP)$
\item[(A')] If $P^\prime(\tau) \in \cal{F}_{i}(\alpha,\beta,OP)$, then $P^\prime(R^2\tau) \in \cal{F}_{i}(\alpha,\alpha+\beta,OP)$
\item[(A')] If $P^\prime(\tau) \in \cal{F}_{i}(\alpha,\beta,RI)$, then $P^\prime(R^2\tau) \in \cal{F}_{i}(\alpha,\alpha+\beta,RI)$
\end{itemize}
\end{lemma}  

\begin{lemma}\label{lem:poch2} Composites of top and right twists act as follows on generating functions in almost quiver form:
\begin{itemize}
\item[(B)] If $P^\prime(\tau) \in \cal{F}_{a}(\alpha,\beta,UP)$, then $P^\prime(RT\tau) \in \cal{F}_{i}(2\alpha+\beta,\alpha+\beta,OP)$
\item[(B)] If $P^\prime(\tau) \in \cal{F}_{a}(\alpha,\beta,OP)$, then $P^\prime(RT\tau) \in \cal{F}_{i}(2\alpha+\beta,\alpha+\beta,RI)$
\item[(B')] If $P^\prime(\tau) \in \cal{F}_{i}(\alpha,\beta,OP)$, then $P^\prime(TR\tau) \in \cal{F}_{a}(\alpha+\beta,\alpha+2\beta,UP)$
\item[(B')] If $P^\prime(\tau) \in \cal{F}_{i}(\alpha,\beta,RI)$, then $P^\prime(TR\tau) \in \cal{F}_{a}(\alpha+\beta,\alpha+2\beta,OP)$
\end{itemize}
\end{lemma}
\begin{lemma}\label{lem:poch3} We also have the following transformations:
\begin{itemize}
\item[(C)] If $P^\prime(\tau) \in \cal{F}_{i}(\alpha,\beta,X)$, then $P^\prime(T\tau) \in \cal{F}_{i}(\alpha+2\beta,\beta,TX)$
\item[(C')] If $P^\prime(\tau) \in \cal{F}_{a}(\alpha,\beta,X)$, then $P^\prime(R\tau) \in \cal{F}_{a}(\alpha,2\alpha+\beta,RX)$
\end{itemize}
where $X$ can be any object type, and $TX$ and $RX$ are its images under those twist operations.
\end{lemma}
We postpone the proofs of these lemmas until Section~\ref{sec:prooflemmas}.

\begin{proof}[Proof of Proposition~\ref{prop:newsumalgo}]
We use the continued fraction expansion of $u/v$ to write $\tau_{u/v}=A_k\ldots A_1 \tau_{0/1}$, where each $A_i$ is either $T$ or $R$, with $A_1=A_k=T$. Then $\tau_{(u-v)/v}=A_{k-1}\ldots A_1\tau_{0/1}$, and we shall study  $P'(\tau_{r_i/s_i})$ for certain intermediate tangles $\tau_{r_i/s_i}:=A_{n_i}\ldots A_1\tau_{0/1}$, where $n_{i+1}=n_i+1$ or $n_i+2$. More specifically, we inductively define intermediate tangles $\tau_{r_i/s_i}$ such that:
\begin{equation}
\label{eqn:condi}
\begin{cases}
P'(\tau_{r_i/s_i})\in \cal{F}_{a}(r_i/2,s_i,X)  \textrm{ with }X\neq RI &\textrm{ if } r_i \textrm{ is even and }s_i\textrm{ is odd (\textit{state T}),} \\
P'(\tau_{r_i/s_i})\in \cal{F}_{i}(r_i,s_i/2,X)\textrm{ with }X \neq UP & \textrm{ if }r_i\textrm{ is odd and }s_i\textrm{ even (\textit{state R})}.
\end{cases}
\end{equation}
 For this, we start with $n_0=0$ and $r_0/s_0=0/1$, for which we evidently have $P'(\tau_{0/1})\in \cal{F}_{a}(0,1,UP)$. Assuming that $\tau_{r_i/s_i}$ satisfies \eqref{eqn:condi}, we define $\tau_{r_{i+1}/s_{i+1}}$ inductively by $n_{i+1}=n_i+2$ if the state of $P'(\tau_{r_i/s_i})$ matches $A_{n_i+1}$, or $n_{i+1}=n_i+1$ if they don't match. In the first case, the relations from Lemmas~\ref{lem:poch1} and \ref{lem:poch2} imply that \eqref{eqn:condi} is satisfied for $P'(\tau_{r_{i+1}/s_{i+1}})$, while in the second case this follows from Lemma~\ref{lem:poch3}.
For example, suppose that $P^\prime(\tau_{r/s})\in \cal{F}_{a}(r/2,s,X)$ with $X \neq RI$, $r$ even and $s$ odd. Then transformation (A) guarantees that $P^\prime(T^2\tau_{r/s})=P^\prime(\tau_{(r+2s)/s})\in \cal{F}_{a}((r+2s)/2,s,X)$ with $X \neq RI$, which is the desired outcome since $r+2s$ is even and $s$ is still odd.  Transformation (B), on the other hand, would produce $P^\prime(RT\tau_{r/s})=P^\prime(\tau_{(r+s)/(r+2s)})\in \cal{F}_{i}(r+s,(r+2s)/2,RTX)$, which also satisfies \eqref{eqn:condi}. In these two examples we have used that $T\tau_{r/s}=\tau_{(r+s)/s}$ and $R\tau_{r/s}=\tau_{r/(r+s)}$.

Now suppose that our inductive construction of intermediate tangles fails to reach $\tau_{(u-v)/v}$. This could only happen if it gets stuck at $A_{k-2}\ldots A_1\tau_{0/1}$, which can either be $\tau_{(u-2v)/v}$ in state $T$ or $\tau_{(u-v)/(2v-u)}$ in state $R$, depending on the type of the remaining crossing $A_{k-1}$, to which none of the transformations in Lemmas~\ref{lem:poch1} -- \ref{lem:poch3} can be applied. However, in such a case, \eqref{eqn:condi} implies that $u$ is even, contradicting the initial assumption of this proposition. Hence, the inductive construction does in fact reach $\tau_{(u-v)/v}$.  

Finally, note that $\tau_{(u-v)/v}$ is of type $UP$ or $RI$ depending on whether $v$ is odd or even since $\tau_{u/v}$ closes to a knot. This leaves only the options $P^\prime(\tau_{(u-v)/v}) \in \cal{F}_{a}((u-v)/2,v,UP)$ or $\cal{F}_{i}(u-v,v/2,RI)$, which proves the proposition.
\end{proof}

\begin{lemma}\label{lemmaclosure} We have:
\[
Cl (T UP[j,k]) 
=(-q)^{k-j} a^{-j} q^{2k^2+j^2} {j \brack k}_+ \frac{ \qp{a^2q^{2-2j-2k}}{k}}{\qpp{k}}
\]
\end{lemma}

\begin{lemma}\label{lemmaclosure2} We have:
\[
Cl_{NS} (T RI[j,k]) 
 =(-q)^k a^{2k-j} q^{-4kj +2k^2+j^2}{j \brack k}_+  \frac{(a^2 q^{2-2j-2(j-k)};q^2)_{j-k}}{(q^2;q^2)_{j-k}}
\]
\end{lemma}
We again postpone the proofs of these lemmas to Section~\ref{sec:prooflemmas}.

\begin{proof}[Proof of Theorem~\ref{thm:knotty}]
Since $K_{u/v}$ is a knot, $u$ is odd. We now have two cases depending on the parity of $v$. For odd $v$, Proposition~\ref{prop:newsumalgo} shows that $P^\prime(\tau_{(u-v)/v})$ is in $\cal{F}_{a}((u-v)/2,v,UP)$ and thus can be written in almost quiver form, with $(u-v)/2$ active and $v$ inactive summation indices and an extra $q$-Pochhammer factor $\qpp{k}$ in each summand, where $k$ denotes the sum of all active summation indices. In order to compute the generating function for the colored HOMFLY-PT polynomials of $K_{u/v}$ we have to account for the final top twist and the closure operation. According to Lemma~\ref{lemmaclosure}, this can be accomplished by replacing each instance of $UP[j,k]$ in $P^\prime(\tau_{(u-v)/v})$ by 
\[(-q)^{k-j} a^{-j} q^{2k^2+j^2}  \frac{ \qp{a^2q^{2-2j-2k}}{k}}{\qpp{k}} x^j,\] since the $q$-binomial ${j \brack k}_+$ is already accounted for in the rescaling of the skein module elements. The denominator $\qpp{k}$ can be cancelled against the extra factor present in $P^\prime(\tau_{(u-v)/v})$, and we obtain an expression that would be in polynomial quiver form, were it not for its additional factors $\qp{a^2q^{2-2j-2k}}{k}$ in the summands. These can be removed by Lemma~\ref{lemlong} at the expense of splitting the $(u-v)/2$ active summation indices into twice as many. The result is an expression for the generating function in polynomial quiver form with $(u-v)+v$ summation indices. More explicitly, the final top twist and the closure act by:
\[   \scalemath{0.75}{ UP, \left( \begin{array}{c|c|c} 1 & S_+ & A_+ \\ \hline\hline 0 & S_- & A_- 
\end{array} \right)},
 \scalemath{0.75}{
 \left( 
\begin{array}{c|c} Q_{++} & Q_{+-} \\ \hline Q_{-+} & Q_{--} 
\end{array} \right)} 
\xrightarrow{} 
\scalemath{0.75}{\left( \begin{array}{c|c} S_+ & A_+-1  \\ \hline  S_--1 & A_--1 \\\hline S_++1 & A_++1 
\end{array} \right)},
 \scalemath{0.75}{
 \left( 
\begin{array}{c|c|c} 
Q_{++}+3 & Q_{+-}+1 &Q_{+-}+1+L \\ \hline
Q_{-+}+1 & Q_{--}+1& Q_{-+}  \\ \hline 
Q_{++}+1+U & Q_{+-}& Q_{++} 
\end{array} \right)} \] 
Here we use a pair notation for the HOMFLY-PT generating function in polynomial quiver form, in analogy with the triple notation for generating functions of skein module elements.

The case of even $v$ is very similar. Proposition~\ref{prop:newsumalgo} shows $P^\prime(\tau_{(u-v)/v}) \in \cal{F}_{i}(u-v,v/2,RI)$ and the last twist and the closure are accomplished by Lemma~\ref{lemmaclosure2}. The factors $\qpp{j-k}$ in the numerator and denominator cancel, and the extra $q$-Pochhammer symbol $(a^2 q^{2-2j-2(j-k)};q^2)_{j-k}$ can be absorbed at the expense of doubling the number of inactive summation indices from $v/2$ to $v$. Here, the final top twist and the closure act by:
\[  \scalemath{0.75}{ OP, \left( \begin{array}{c|c|c} 0 & S_+ & A_+ \\ \hline\hline 1 & S_- & A_- 
\end{array} \right)},
 \scalemath{0.75}{
 \left( 
\begin{array}{c|c} Q_{++} & Q_{+-} \\ \hline Q_{-+} & Q_{--} 
\end{array} \right)} \xrightarrow{}
 \scalemath{0.75}{\left( \begin{array}{c|c} S_++1 & A_++1 \\\hline S_- & A_--1 \\ \hline  S_-+1 & A_-+1 
\end{array} \right)},
 \scalemath{0.75}{
 \left( 
\begin{array}{c|c|c} 
Q_{++}-1 &Q_{+-}-1 & Q_{+-}-2 \\ \hline 
Q_{-+}-1 & Q_{--}+1& Q_{--}-1+L \\ \hline
Q_{-+}-2 & Q_{--}-1+U & Q_{--}-2
\end{array} \right)} \] 
The final result is again in polynomial quiver form with $u$ summation indices.
\end{proof}

\begin{proof}[Proof of Theorem~\ref{thm:linky}]
We follow the algorithm in Theorem~\ref{thm:links-full} to bring $P^\prime(\tau_{(u-v)/v})$ into quiver form. A simple inductive argument shows that the result has $u-v$ active and $v$ inactive summation indices. Then we use Lemma~\ref{lemmaclosure} or Lemma~\ref{lemmaclosure2} to perform the last top twist and the closure. This produces a generating function with $u$ summation indices, which would be in quiver form, were it not for additional factors
\[\qpp{j}\frac{\qp{a^2q^{2-2j-2k}}{k}}{\qpp{k}}=\qp{q^{2+2k}}{j-k} \qp{a^2q^{2-2j-2k}}{k} \]
or \[ \qpp{j}\frac{(a^2 q^{2-2j-2(j-k)};q^2)_{j-k}}{(q^2;q^2)_{j-k}}=\qp{q^{2+2j-2k}}{k}(a^2 q^{2-2j-2(j-k)};q^2)_{j-k}\]
in the summands. In both cases, these factors can be absorbed using Lemma~\ref{lemlong}, which doubles the number of summation indices from $u$ to $2u$.
\end{proof}

For knots $K_{u/v}$, Theorem~\ref{thm:links} can be deduced from Theorem~\ref{thm:knots}. Indeed, the generating function in polynomial quiver form \eqref{eqn:knottyquiverform} with $u$ summation indices, provided by the latter, transforms into a generating function in quiver form \eqref{eqn:quiverform} with $2u$ summation indices upon removing the extra $q$-Pochhammer symbol using Lemma~\ref{lemlong}.

On the other hand, generating functions for 2-component links $L_{u/v}$ cannot be brought into polynomial quiver form, simply because the reduced colored HOMFLY-PT polynomials are not Laurent polynomials, see Remark~\ref{rem:HOMFLY-links}.

\subsection{Proofs of Lemmas}
\label{sec:prooflemmas}
In this section we will supply the proofs of the lemmas that were stated and used in the previous section. In these proofs we will use the following lemmas to reduce the number of summation indices of a generating function in almost quiver form, at the expense of creating a new $q$-Pochhammer symbol in the numerator. 
\begin{lemma}\label{lem:csumpre} For $e\geq 0$ we have:
\[\sum_{b,c\geq 0} \frac{(-q)^b q^{b^2+2 e b}}{\qpp{b}\qpp{c}} x^{b+c} = \sum_{d\geq 0} \frac{\qp{q^{2+2e}}{d}}{\qpp{d}} x^d = \sum_{d\geq 0} \frac{\qpp{d+e}}{\qpp{d}\qpp{e}} x^d\]
\end{lemma}
\begin{proof}
This is a special case of Lemma~\ref{lemlong}.\end{proof}

\begin{lemma}
\label{lem:csum} The following represents an equality of generating functions:
\[
\scalemath{0.75}{ X, \left( \begin{array}{c|c|c} 1 & S_1 & A_1 \\ \hline 0 & S_2+1 & A_2\\ \hline 0 & S_2 & A_2  \\ \hline 0 & S_3 & A_3 
\end{array} \right)},
 \scalemath{0.75}{
 \left( 
\begin{array}{c|c|c|c}
 Q_{11} & Q_{12}+1& Q_{12}& Q_{13} \\ \hline 
 Q_{21}+1 & Q_{22}+1& Q_{22}+L& Q_{23} \\ \hline 
 Q_{21} & Q_{22}+U& Q_{22}& Q_{23} \\ \hline 
 Q_{31} & Q_{32}& Q_{32}& Q_{33} 
\end{array} \right)} 
\quad = \quad
\scalemath{0.75}{ X, \left( \begin{array}{c|c|c} 1 & S_1 & A_1 \\ \hline 1 & S_2 & A_2 \\ \hline 0 & S_3 & A_3
\end{array} \right)},
 \scalemath{0.75}{
 \left( 
\begin{array}{c|c|c}
 Q_{11} & Q_{12}& Q_{13} \\ \hline 
 Q_{21} & Q_{22}& Q_{23} \\ \hline 
 Q_{31} & Q_{32}& Q_{33} 
\end{array} \right)}
\]
Here we assume that the second and the third block of summation indices on the left-hand side are both active or both inactive. The other blocks may be empty.
\end{lemma}
\begin{proof}
It suffices to prove this in the case of when the size of the blocks to be summed is one. The general case follows by induction. We denote the summation indices corresponding to the second and third block by $b$ and $c$, the vectors of summation indices in the first and last block by $\textbf{e}$ and $\textbf{f}$, and their sums by $\underline{e}$ and $\underline{f}$ respectively. The relevant part of the generating function on the left-hand side that depends on $b$ and $c$ can be isolated as:
\[\qpp{\underline{e}}\sum_{b,c\geq 0} (-q)^{S_2(b+c)} a^{A_2(b+c)} q^{Q_{22}(b+c)^2+2 (b+c)Q_{12}\cdot\textbf{e}+2 (b+c)Q_{23}\cdot \textbf{f}}\frac{(-q)^b q^{b^2+2 \underline{e} b}}{\qpp{b}\qpp{c}}\]
Using Lemma~\ref{lem:csumpre}, this can be expressed in terms of the single summation index $d=b+c$.
\[\qpp{\underline{e}}\sum_{d\geq 0} (-q)^{S_2 d} a^{A_2 d} q^{Q_{22}d^2+2 d Q_{12}\cdot \textbf{e}+2 d Q_{23}\cdot \textbf{f}}\frac{\qpp{d+\underline{e}}}{\qpp{d}\qpp{\underline{e}}}\]
This is precisely the relevant part of the generating function on the right-hand side.

\end{proof}

\begin{proof}[Proof of Lemma~\ref{lem:poch1}]
We first treat the case of (A), starting with type $UP$, for which we assume the presence of an additional $q$-Pochhammer symbol $\qpp{k}$ where $k$ is the sum of the active summation indices. To this end, let the initial generating function be encoded by the following triple: 
\[  \scalemath{0.75}{ UP, \left( \begin{array}{c|c|c} 1 & S_+ & A_+ \\ \hline\hline 0 & S_- & A_- 
\end{array} \right)},
 \scalemath{0.75}{
 \left( 
\begin{array}{c|c} Q_{++} & Q_{+-} \\ \hline Q_{-+} & Q_{--} 
\end{array} \right)} \] 
Interpreting this as a generating function for a rescaled skein module element, we can apply the twist rule $T$ as in Proposition~\ref{prop:twisttriple} twice to obtain the following:

\begin{equation}\label{eqn:beforesum}  \scalemath{0.75}{ UP, \left( \begin{array}{c|c|c} 1 & S_+ & A_+ \\ \hline 0 & S_- & A_-\\ \hline 0 & S_- -1 & A_-  \\ \hline\hline 0 & S_- -2 & A_- 
\end{array} \right)},
 \scalemath{0.75}{
 \left( 
\begin{array}{c|c|c|c}
 Q_{++}+2 & Q_{+-}+1& Q_{+-}& Q_{+-} \\ \hline 
 Q_{-+}+1 & Q_{--}+1& Q_{--}+L& Q_{--}+L \\ \hline 
 Q_{-+} & Q_{--}+U& Q_{--}& Q_{--}+L \\ \hline 
 Q_{-+} & Q_{--}+U& Q_{--}+U& Q_{--} 
\end{array} \right)} \end{equation}
The second and third block can now be summed using Lemma~\ref{lem:csum}, and we obtain:
 \[  \scalemath{0.75}{ UP, \left( \begin{array}{c|c|c} 1 & S_+ & A_+ \\ \hline 1 & S_- -1 & A_-  \\ \hline\hline 0 & S_- -2 & A_- 
\end{array} \right)},
 \scalemath{0.75}{
 \left( 
\begin{array}{c|c|c}
 Q_{++}+2 & Q_{+-}& Q_{+-} \\ \hline 
 Q_{-+} & Q_{--}& Q_{--}+L \\ \hline 
 Q_{-+} & Q_{--}+U& Q_{--} 
\end{array} \right)} \]

Thus we have proven that $T^2$ transforms generating functions of type $UP$ as claimed in Lemma~\ref{lem:poch1}. 

The case of $OP$ is completely analogous and results in:
 \[  \scalemath{0.75}{ OP, \left( \begin{array}{c|c|c} 1 & S_+ +2 & A_+ +2 \\ \hline 1 & S_- +1 & A_-+1  \\ \hline\hline 0 & S_-  & A_- 
\end{array} \right)},
 \scalemath{0.75}{
 \left( 
\begin{array}{c|c|c}
 Q_{++}-2 & Q_{+-}-3& Q_{+-}-2 \\ \hline 
 Q_{-+}-3 & Q_{--}-2& Q_{--}-1+L \\ \hline 
 Q_{-+}-2 & Q_{--}-1+U& Q_{--} 
\end{array} \right)} \]
The transformations (A') for $OP$ and $RI$ are also similar. The main difference is that we start with:

\[  \scalemath{0.75}{ X, \left( \begin{array}{c|c|c} 0 & S_+ & A_+ \\ \hline\hline 1 & S_- & A_- 
\end{array} \right)},
 \scalemath{0.75}{
 \left( 
\begin{array}{c|c} Q_{++} & Q_{+-} \\ \hline Q_{-+} & Q_{--} 
\end{array} \right)} \]  
The end result for $OP$ is:
 \[  \scalemath{0.75}{ OP, \left( \begin{array}{c|c|c} 0 & S_+ & A_+ \\ \hline\hline 1 & S_+ -2 & A_+-1  \\ \hline 1 & S_--2  & A_--2 
\end{array} \right)},
 \scalemath{0.75}{
 \left( 
\begin{array}{c|c|c}
 Q_{++} & Q_{++}+1+L& Q_{+-}+2 \\ \hline 
 Q_{++}+1+U & Q_{++}+1& Q_{+-}+1 \\ \hline 
 Q_{-+}+2 & Q_{-+}+1& Q_{--}+2 
\end{array} \right)} \]
The end result for $RI$ is:
 \[  \scalemath{0.75}{ RI, \left( \begin{array}{c|c|c} 0 & S_++2 & A_+ \\ \hline\hline 1 & S_+  & A_+  \\ \hline 1 & S_-  & A_- 
\end{array} \right)},
 \scalemath{0.75}{
 \left( 
\begin{array}{c|c|c}
 Q_{++} & Q_{++}+L& Q_{+-} \\ \hline 
 Q_{++}+U & Q_{++}-1& Q_{+-}-2 \\ \hline 
 Q_{-+} & Q_{-+}-2& Q_{--}-2 
\end{array} \right)} \]
These resulting generating functions are of the form claimed in Lemma~\ref{lem:poch1}. 
\end{proof}

\begin{proof}[Proof of Lemma~\ref{lem:poch2}]
The proof is similar to the last one, but slightly more involved. We first treat the transformation (B') in the case of $OP$. To this end, let the initial generating function be encoded by the following triple: 
\[  \scalemath{0.75}{ OP, \left( \begin{array}{c|c|c} 0 & S_+ & A_+ \\ \hline\hline 1 & S_- & A_- 
\end{array} \right)},
 \scalemath{0.75}{
 \left( 
\begin{array}{c|c} Q_{++} & Q_{+-} \\ \hline Q_{-+} & Q_{--} 
\end{array} \right)} \] 

After applying a right twist $ROP$ and then a top twist $TUP$ we get:
\begin{equation} \scalemath{0.75}{
\left( \begin{array}{c|c|c} 
0& S_+&A_+\\ \hline
0& S_+-1&A_+\\ \hline
1& S_--1&A_--1\\ \hline \hline
0& S_+-2&A_+\\ \hline
1& S_--2&A_--1
\end{array} \right)}
,
 \scalemath{0.75}{
 \left( 
\begin{array}{c|c|c|c|c} 
Q_{++}+1 & Q_{++}+L & Q_{+-}+1& Q_{++}+L & Q_{+-}+1 \\ \hline 
Q_{++} + U& Q_{++} & Q_{+-}& Q_{++}+L & Q_{+-}  \\ \hline 
Q_{-+}+1 & Q_{-+} & Q_{--}+1 & Q_{-+}+1 & Q_{--}+1+L \\ \hline 
Q_{++} + U& Q_{++}+U & Q_{+-}+1& Q_{++} & Q_{+-}  \\ \hline 
Q_{-+}+1 & Q_{-+} & Q_{--}+1+U & Q_{-+} & Q_{--}+1
\end{array} \right)}
\end{equation} 
 The indices in the first two blocks can now be summed up in pairs using Lemma~\ref{lem:csum}, which yields:
\[ \scalemath{0.75}{
\left( \begin{array}{c|c|c} 
1& S_+-1&A_+\\ \hline
1& S_--1&A_--1\\ \hline \hline
0& S_+-2&A_+\\ \hline
1& S_--2&A_--1
\end{array} \right)}
, 
 \scalemath{0.75}{
 \left( 
\begin{array}{c|c|c|c} 
 Q_{++} & Q_{+-}& Q_{++}+L & Q_{+-}  \\ \hline 
 Q_{-+} & Q_{--}+1 & Q_{-+}+1 & Q_{--}+1+L \\ \hline 
 Q_{++}+U & Q_{+-}+1& Q_{++} & Q_{+-}  \\ \hline 
 Q_{-+} & Q_{--}+1+U & Q_{-+} & Q_{--}+1
\end{array} \right)
}
 \] 

The new $K$ vector $(1|1|0|1)$ corresponds to an additional $q$-Pochhammer symbol $\qpp{h+x}$ where $h$ is the new active index given by the sum of the first two blocks of summation indices and $x$ is the sum of the indices in the last block. This additional factor is larger than we want, so we split off the extra term $\qpp{h+x}= \qpp{h}\qp{q^{2+2h}}{x}$. This extra term can be absorbed into the quiver sum at the expense of splitting the summation indices in the fourth block, see Lemma~\ref{lemlong}. So, finally, we get:
\[ \scalemath{0.75}{
\left( \begin{array}{c|c|c} 
1& S_+-1&A_+\\ \hline
1& S_--1&A_--1\\ \hline \hline
0& S_+-2&A_+\\ \hline
0& S_--2&A_--1\\ \hline
0& S_--1&A_--1
\end{array} \right)}
,
 \scalemath{0.75}{
 \left( 
\begin{array}{c|c|c|c|c} 
 Q_{++} & Q_{+-}& Q_{++}+L & Q_{+-}& Q_{+-}+1  \\ \hline 
 Q_{-+} & Q_{--}+1 & Q_{-+}+1 & Q_{--}+1+L& Q_{--}+2+L \\ \hline 
 Q_{++}+U & Q_{+-}+1& Q_{++} & Q_{+-}& Q_{+-}  \\ \hline 
 Q_{-+} & Q_{--}+1+U & Q_{-+} & Q_{--}+1& Q_{--}+1+U\\ \hline 
 Q_{-+}+1 & Q_{--}+2+U & Q_{-+} & Q_{--}+1+L& Q_{--}+2
\end{array} \right)}
 \] 
This is of the desired form. Had we started with $RI$ instead, we would have ended up with the following:
\[ \scalemath{0.75}{
\left( \begin{array}{c|c|c} 
1& S_+ + 1&A_+ + 1\\ \hline
1& S_- + 1&A_- + 1\\ \hline \hline
0& S_+&A_+\\ \hline
0& S_-&A_-\\ \hline
0& S_- + 1&A_-
\end{array} \right)}
, 
 \scalemath{0.75}{
 \left( 
\begin{array}{c|c|c|c|c} 
 Q_{++}-2 & Q_{+-}-3& Q_{++}-1+L & Q_{+-}-2& Q_{+-}-1  \\ \hline 
 Q_{-+}-3 & Q_{--}-3 & Q_{-+}-1 & Q_{--}-2+L& Q_{--}-1+L \\ \hline 
 Q_{++}-1+U & Q_{+-}-1& Q_{++} & Q_{+-}-1& Q_{+-}-1  \\ \hline 
 Q_{-+}-2 & Q_{--}-2+U & Q_{-+}-1 & Q_{--}-1& Q_{--}-1+U\\ \hline 
 Q_{-+}-1 & Q_{--}-1+U & Q_{-+}-1 & Q_{--}-1+L& Q_{--}
\end{array} \right)}
 \] 
Next, we treat the transformation (B) in the case $UP$, where we start with:
\[  \scalemath{0.75}{ \left( \begin{array}{c|c|c} 1 & S_+ & A_+ \\ \hline\hline 0 & S_- & A_- 
\end{array} \right)}
,  \scalemath{0.75}{\left( 
\begin{array}{c|c} Q_{++} & Q_{+-} \\ \hline Q_{-+} & Q_{--} 
\end{array} \right)} \]  
After a top twist $TUP$ and then a right twist $RUP$ we get: 
\begin{equation}
\label{eqn:stepp} \scalemath{0.75}{\left( \begin{array}{c|c|c} 
1& S_+ & A_+ \\ \hline
0& S_- & A_- \\ \hline\hline
1& S_+-1 & A_+-1 \\ \hline
0& S_--1 & A_- -1 \\ \hline
0 & S_- - 2 &A_- -1
\end{array} \right)}
,
 \scalemath{0.75}{
 \left( 
\begin{array}{c|c|c|c|c} 
Q_{++} +1 & Q_{+-} & Q_{++} +2+L & Q_{+-}+1 & Q_{+-}+1 \\ \hline 
Q_{-+} & Q_{--} & Q_{-+}+2 & Q_{--}+1+L & Q_{--}+1+L \\ \hline 
Q_{++} +2+U & Q_{+-}+2 & Q_{++} +3 & Q_{+-}+2 & Q_{+-}+1 \\ \hline 
Q_{-+}+1 & Q_{--}+1+U & Q_{-+}+2 & Q_{--}+2 & Q_{--}+1+L \\ \hline 
Q_{-+}+1 & Q_{--}+1+U & Q_{-+}+1 & Q_{--}+1+U & Q_{--}+1 
\end{array} \right)} 
\end{equation}
We would like to use Lemma~\ref{lem:csum} to sum the indices in the last two blocks in pairs to get:
\begin{equation}
\label{eqn:step} \scalemath{0.75}{ \left( \begin{array}{c|c|c} 
0& S_+ & A_+ \\ \hline
0& S_- & A_- \\ \hline\hline
1& S_+-1 & A_+-1 \\ \hline
1& S_--2 & A_- -1
\end{array} \right)}
,
 \scalemath{0.75}{
 \left( 
\begin{array}{c|c|c|c} 
Q_{++} +1 & Q_{+-} & Q_{++} +2+L & Q_{+-}+1  \\ \hline 
Q_{-+} & Q_{--} & Q_{-+}+2 & Q_{--}+1+L  \\ \hline 
Q_{++} +2+U & Q_{+-}+2 & Q_{++} +3 & Q_{+-}+1 \\ \hline 
Q_{-+}+1 & Q_{--}+1+U & Q_{-+}+1 & Q_{--}+1 
\end{array} \right)} 
\end{equation}
However, this would require that the column vector $K$ in \eqref{eqn:stepp} be $(0|0|1|0|0)$ instead of $(1|0|1|0|0)$. Thus we instead obtain \eqref{eqn:step} multiplied by a quotient of $q$-Pochhammer symbols of type $(1|0|1|0) $ and $(0|0|1|0)$, which will split the summation indices in the first block via Lemma~\ref{lemlong} to give:
\[ \scalemath{0.75}{ \left( \begin{array}{c|c|c} 
0& S_+ & A_+ \\ \hline
0& S_++1 & A_+ \\ \hline
0& S_- & A_- \\ \hline\hline
1& S_+-1 & A_+-1 \\ \hline
1& S_--2 & A_- -1
\end{array} \right)
}
,
 \scalemath{0.75}{
 \left( 
\begin{array}{c|c|c|c|c} 
Q_{++} +1 & Q_{++} +1+U & Q_{+-} & Q_{++} +2+L & Q_{+-}+1  \\ \hline 
 Q_{++} +1+L &Q_{++} +2 & Q_{+-} & Q_{++} +3+L & Q_{+-}+1  \\ \hline 
Q_{-+} &Q_{-+} & Q_{--} & Q_{-+}+2 & Q_{--}+1+L  \\ \hline 
Q_{++} +2+U &Q_{++} +3+U & Q_{+-}+2 & Q_{++} +3 & Q_{+-}+1 \\ \hline 
Q_{-+}+1 &Q_{-+}+1 & Q_{--}+1+U & Q_{-+}+1 & Q_{--}+1 
\end{array} \right)
} \] 
which is of the desired form. The case of $OP$ is analogous, with end result:
\[ 
 \scalemath{0.75}{
\left( \begin{array}{c|c|c} 
0& S_+ +2 & A_+ +1 \\ \hline
0& S_++3 & A_+ +1 \\ \hline
0& S_-+2 & A_- \\ \hline\hline
1& S_+ +1& A_+ +1 \\ \hline
1& S_- & A_- 
\end{array} \right)
}
,
 \scalemath{0.75}{
 \left( 
\begin{array}{c|c|c|c|c} 
Q_{++} -1 & Q_{++} -1+U & Q_{+-}-1 & Q_{++} -1+L & Q_{+-}-1 \\ \hline 
 Q_{++} -1+L &Q_{++}  & Q_{+-}-1 & Q_{++} +L & Q_{+-}-1 \\ \hline 
Q_{-+}-1 &Q_{-+}-1 & Q_{--} & Q_{-+} & Q_{--}+L  \\ \hline 
Q_{++} -1+U &Q_{++} +U & Q_{+-} & Q_{++} -1 & Q_{+-}-2 \\ \hline 
Q_{-+}-1 &Q_{-+}-1 & Q_{--}+U & Q_{-+}-2 & Q_{--}-1 
\end{array} \right) }\] 
\end{proof}

We have included the end results of all steps in the last two proofs, because they can be used to efficiently code an algorithm that computes the desired generating functions. Conversely, the steps in the proof can be quickly checked using a computer implementation of the twist rules, Lemma~\ref{lem:csum} and Lemma~\ref{lemlong}. 

\begin{proof}[Proof of Lemma~\ref{lem:poch3}]
We first use the appropriate twist rule from Proposition~\ref{prop:twisttriple}. The resulting expressions are in almost quiver form, however, the extra $q$-Pochhammer symbols are too long. For example, the result of a top twist applied to $P^\prime(\tau) \in \cal{F}_{i}(\alpha,\beta,X)$ would have $\alpha+\beta$ active and $\beta$ inactive summation indices, but also carry extra factors $\qpp{j-k}$, whose length $j-k$ is the sum of the old inactive summation indices. The desired extra factor is $\qpp{j-h}$, whose length is the sum of the new inactive summation indices. Their quotient is
$\qpp{j-k}/\qpp{j-h}=\qp{q^{2+2j-2h}}{h-k}$, which doubles the number of new active summation indices to $2\beta$ via Lemma~\ref{lemlong}. In fact, this is an instance of Lemma~\ref{lem:csum} read from right to left, so we omit further details. The other five cases are analogous.
\end{proof}

\begin{proof}[Proof of Lemma~\ref{lemmaclosure}]
We start by expanding $Cl (T UP[j,k])$:
\begin{align*}
Cl (T UP[j,k]) &= \sum_{h=k}^j (-q)^{h-j} q^{k^2} {h \brack k}_+ a^{-j} q^{j^2+h^2} 
\frac{\qp{a^2 q^{2-2j-2h}}{j}}{\qpp{j}} {j\brack h}_+
\end{align*}
Comparing with the right-hand side of the proposed identity, and using ${h \brack k}_+{j\brack h}_+={j\brack k}_+{j-k\brack h-k}_+$, we find that it remains to prove:
\begin{equation}
\label{eqn:1}\sum_{h=k}^j (-q)^{h} q^{h^2} {j-k \brack h-k}_+ \frac{ \qp{a^2q^{2-2j-2h}}{j}}{\qpp{j}} = (-q)^{k} q^{k^2} \frac{ \qp{a^2q^{2-2j-2k}}{k}}{\qpp{k}} 
\end{equation}
Here we rewrite the left-hand side in terms of the variable $x=h-k$:
\begin{gather*}
\sum_{x=0}^{j-k} (-q)^{k+x} q^{k^2 + 2 kx + x^2} {j-k \brack x}_+ \frac{ \qp{a^2q^{2-2j-2k-2x}}{j}}{\qpp{j}}.
\end{gather*}
Now we use Lemma~\ref{mainlemma}:
\begin{gather}
= \sum_{x=0}^{j-k}  \frac{(-q)^{k+x} q^{k^2 + 2 kx + x^2}}{\qpp{j}}{j-k \brack x}_+  \sum_{y=0}^{j} (-1)^y a^{2y} q^{(2-2j-2k-2x)y +y^2-y} {j \brack y}_+ \nonumber \\
 =\frac{(-q)^{k} q^{k^2}}{\qpp{j}} \sum_{y=0}^{j} (-q)^y a^{2y} q^{(-2j-2k)y +y^2} {j \brack y}_+ \sum_{x=0}^{j-k}  (-1)^{x} q^{2x (1+k-y) +x^2-x }{j-k \brack x}_+ 
\label{eqn:stst}
\end{gather}
Now we use Lemma~\ref{mainlemma} again:
\[\sum_{x=0}^{j-k}  (-1)^{x} q^{2x (1+k-y) +x^2-x }{j-k \brack x}_+ = \qp{q^{2(1+k-y)}}{j-k} =\]
\[= \delta_{\{y\leq k\}} \frac{\qpp{j-y}}{\qpp{k-y}}.\]
Plugging this back into \eqref{eqn:stst} gives:
\begin{align*}
& \frac{(-q)^{k} q^{k^2}}{\qpp{j}} \sum_{y=0}^{k} (-1)^y a^{2y} q^{(2-2j-2k)y +y^2-y} \frac{\qpp{j}}{\qpp{y} \qpp{k-y}}=\\
=& \frac{(-q)^{k} q^{k^2}}{\qpp{k}} \sum_{y=0}^{k} (-1)^y a^{2y} q^{(2-2j-2k)y +y^2-y}  {k \brack y}_+.
\end{align*}
One final application of Lemma~\ref{mainlemma} verifies \eqref{eqn:1} and completes the proof.
\end{proof}

\begin{proof}[Proof of Lemma~\ref{lemmaclosure2}]
We start by expanding $Cl_{NS} (T RI[j,k])$:
\begin{align*}
Cl_{NS} (T RI[j,k])&=\sum_{h=k}^j (-q)^h a^h  q^{k^2-2jh} {h \brack k}_+ a^{-j+h} q^{(j-h)^2} {j \brack h}_{+} \frac{(a^2 q^{2-2j};q^2)_{j-h}}{(q^2;q^2)_{j-h}}
\end{align*}
Let $x=h-k$, and define:
$$L:=\sum_{x=0}^{j-k} (-q)^x a^{2x} q^{-4jx+2kx+x^2} {j-k \brack x}_{+} (a^2 q^{2-2j};q^2)_{j-k-x} (q^{2+2j-2k-2x};q^2)_x$$
Then, since ${j \brack h}_+{h\brack k}_+={j\brack k}_+{j-k \brack x}_+$ the statement of the lemma reduces to proving 
that $L=(a^2q^{2-4j+2k};q^2)_{j-k}$. To this end, first we use the $q$-binomial identities on the two $q$-Pochhammer symbols from $L$:
\begin{eqnarray*}
(a^2 q^{2-2j};q^2)_{j-k-x}&=&\sum_{y=0}^{j-k-x} (-q)^y a^{2y} q^{y^2-2jy} {j-k-x \brack y}_{+},\\
(q^{2+2j-2k-2x};q^2)_x&=&\sum_{z=0}^x (-q)^z q^{z^2+2jz-2kz-2xz} {x\brack z}_{+}.
\end{eqnarray*}
Now we introduce change of variables: $s=x+y$, $\alpha=x-z$ and consider $s$, $\alpha$ and $z$ as new summation indices. Since
${j-k \brack x}_{+} {j-k-x \brack y}_{+} {x\brack z}_{+} ={j-k\brack s}_{+} {s \brack \alpha}_{+} {s-\alpha \brack z}_{+}$, we get:

\begin{align*}L=&\sum_{s=0}^{j-k} (-q)^s a^{2s} q^{s^2-2js} {j-k \brack s}_{+}\sum_{\alpha=0}^s q^{-2j\alpha +2k \alpha-2s \alpha+2 \alpha^2}{s\brack \alpha}_{+} \\
&\qquad \times 
\sum_{z=0}^{s-\alpha} (-q)^z q^{z^2+2\alpha z-2sz} {s-\alpha \brack z}_{+}.
\end{align*}
The innermost summation over $z$ is equal to $(q^{2+2\alpha-2s};q^2)_{s-\alpha}$ by the $q$-binomial identity, and is furthermore equal to $\delta_{\alpha,s}$, since by definition it is zero as soon as $s>\alpha$. Therefore:
$$L=\sum_{s=0}^{j-k} (-1)^s a^{2s} q^{s^2+s-2js} {j-k \brack s}_{+}\sum_{\alpha=0}^s q^{-2j\alpha +2k \alpha-2s \alpha+2 \alpha^2} {s\brack \alpha}_{+}\delta_{\alpha,s}$$
$$\quad=\sum_{s=0}^{j-k} (-1)^s a^{2s} q^{s^2+s-2js} {j-k \brack s}_{+} q^{-2js +2ks}=(a^2q^{2-4j+2k};q^2)_{j-k},$$
again by the $q$-binomial identity, which finishes the proof.
\end{proof}

\subsection{Quivers and HOMFLY-PT homology}
Theorem~\ref{thm:knotsHHH} claims that the reduced HOMFLY-PT homology of a rational knot $K$ can be recovered from the generating function of the colored HOMFLY-PT polynomials obtained in Theorem~\ref{thm:knotty}. This is unsurprising since Rasmussen has shown that the reduced HOMFLY-PT homology of a rational knot is determined by its Euler characteristic and the signature of the knot. In our grading conventions, his result is the following.

\begin{lemma}[{\cite[Corollary 1]{Ras}}] Let $K$ be a rational knot, then the $a,q,t$-trigradings of the generators of the reduced HOMFLY-PT homology of $K$ satisfy $2\mathrm{gr}_t-2\mathrm{gr}_a-\mathrm{gr}_q=\sigma(K)$ where $\sigma(K)$ denotes the signature of $K$, normalized so that positive knots have negative signature.  
\end{lemma}
Since the homological degrees of the generators are determined by the other degrees and the signature, there can be no cancellation in the bigraded Euler characteristic, and so the Poincar\'e polynomial of the HOMFLY-PT homology can be recovered from the reduced HOMFLY-PT polynomial.

In order to prove Theorem~\ref{thm:knotsHHH}, it thus suffices to check that for the generating function computed in Theorem~\ref{thm:knotty} for the rational knot $K$ we have
\[-2 Q_{i,i} - 2 a_i -(-Q_{i,i}-q_i)=q_i-Q_{i,i}-2 a_i = \sigma(K) \text{ for all }i\]
(Note that we implicitly use the $q \leftrightarrow q^{-1}$ symmetry of the reduced uncolored HOMFLY-PT polynomial to identify the decategorification of the HOMFLY-PT homology with the coefficient of $x$ in \eqref{eqn:knotty}, see also Section~\ref{sec:framing}.
For such a generating function, we define the $\delta$-grading of a summation index $i$ to be $\delta_i:=q_i-Q_{i,i}-2 a_i$. Note that this notion also makes sense for the generating functions associated to tangles. We will now show that generating functions from Theorem~\ref{thm:knotty} are homogeneous with respect to the $\delta$-grading. This will be done by induction of the number of crossings on the link, following the building process of the rational tangle. 

\begin{lemma}\label{lem:hom} Let $c$ and $d$ be sums of disjoint sets of summation indices. Then multiplying by a Pochhammer symbol $\qp{a^x q^{2+2cy+2dz}}{d}$ and absorbing it via Lemma~\ref{lemlong} preserves the $\delta$-grading of a $\delta$-homogeneous generating functions in (almost) quiver form provided that $x=-z$.
\end{lemma}
\begin{proof}
Suppose we have a generating function that is $\delta$-homogeneous and that $d=d_1+\cdots d_k$. Then multiplication by the Pochhammer symbol and absorption via Lemma~\ref{lemlong} will split the summation indices $d_i=\alpha_i+\beta_i$ and we get:

\begin{align*}\frac{\qp{a^x q^{2+2cy+2dz}}{d}}{(q^2;q^2)_{d_1}\cdots (q^2;q^2)_{d_k}}
&\\ =\sum\limits_{\substack{\alpha_1+\beta_1=d_1\\\cdots\\\alpha_k+\beta_k=d_k}}&
\frac{(-q a^x q^{2cy+2dz})^{{\alpha_1+\ldots+\alpha_k}} q^{\alpha_1^2+\ldots+\alpha_k^2+2\sum_{i=1}^{k-1} \alpha_{i+1} (d_1+\ldots+d_i)}}{(q^2;q^2)_{\alpha_1}\cdots(q^2;q^2)_{\alpha_k}(q^2;q^2)_{\beta_1}\cdots(q^2;q^2)_{\beta_k}} 
\end{align*} 
Here the linear $q$-degree increases only for the summation indices $\alpha_i$, but this is compensated by the increase in quadratic $q$-degree. The quadratic terms in the $q$-exponent that mix $c$ and $\alpha_i$ as well as the sum on the right only influence off-diagonal entries of the adjacency matrix of the quiver and hence do not contribute to any $\delta_i$. Thus the net change in $\delta$ is $-2z-2x$ for the $\alpha_i$ and zero for the $\beta_i$.
 \end{proof}
 
\begin{corollary}\label{cor:twisthom} The twist operations as detailed in Theorem~\ref{thm:prodform} preserve the $\delta$-homogeneity of the generating functions of rescaled skein module elements. More specifically, $T^\prime OP$, $T^\prime RI$, $R^\prime UP$ and $R^\prime OP$ preserve the $\delta$-grading, whereas $T^\prime UP$ decreases it by one and $R^\prime RI$ increases it by one.
\end{corollary}
\begin{proof} For each of the six twist rules, one easily checks via Lemma~\ref{lem:hom} that multiplication by the $q$-Pochhammer symbol preserves the $\delta$-grading, whereas the monomial factor preserves or shifts the grading as stated.
\end{proof}
The next corollary follows analogously.
\begin{corollary}\label{cor:closehom} The closure operation coupled to the last crossing, as simplified in Lemma~\ref{lemmaclosure} and Lemma~\ref{lemmaclosure2} (and rescaled to be compatible with the algorithm from Theorem~\ref{thm:prodform}), preserves $\delta$-homogeneity. More specifically $Cl (T UP)$ preserves the $\delta$-grading, whereas  $Cl_{NS} (T RI)$ increases it by one. 
\end{corollary}

The following is immediate upon inspection.
\begin{lemma}\label{lem:sumhom} The transformation described in Lemma~\ref{lem:csum} preserves $\delta$-homogeneity.
\end{lemma}

\begin{proposition}
The generating functions for the colored HOMFLY-PT polynomials of rational knots $K_{u/v}$, brought into polynomial quiver form as in Theorem~\ref{thm:knotty}, are $\delta$-homogeneous with $\delta$-grading equal to the signature $\sigma(K_{u/v})$.
\end{proposition}
\begin{proof} We follow a variation of the algorithm described in the proof of Proposition~\ref{prop:newsumalgo}, where we use the rescaled twist operations from Theorem~\ref{thm:prodform} instead of the usual ones. The generating function of the trivial tangle is $\delta$-homogeneous and by Corollary~\ref{cor:twisthom}, this is preserved by the twist rules. In the proof of Proposition~\ref{prop:newsumalgo} we apply these twist rules in pairs and rewrite the resulting generating functions by Lemma~\ref{lem:csum} or by absorbing Pochhammer symbols via Lemma~\ref{lemlong}. These types of rewriting preserve $\delta$-homogeneity by Lemma~\ref{lem:sumhom} and Lemma~\ref{lem:hom} respectively. Finally, the closure operation is homogeneous by Corollary~\ref{cor:closehom}. The final $\delta$-grading is computed as $1+\#\{R^\prime X\} - \#\{T^\prime UP, R^\prime UP,R^\prime OP\}$, where the second term indicates the number of right twists in the rational tangle and the third term the number of occurrences of twists of the form $T^\prime UP$, $R^\prime UP$ and $R^\prime OP$. To identify the $\delta$-grading as the signature, we use the Gordon--Litherland formula
\[\sigma(K_{u/v}) =  \sigma(D_{u/v})-\mu(D_{u/v})\]
for a standard chequerboard-colored diagram $D_{u/v}$ of $K_{u/v}$. Here, $\sigma(D_{u/v})$ denotes the signature of the Goeritz matrix of $D_{u/v}$, which is equal to $1+\#\{R^\prime X\}$ by \cite{QAQ}, and $\mu(D_{u/v})$ is a correction term that is easily seen to agree with $\#\{T^\prime UP, R^\prime UP,R^\prime OP\}$ from the original description in \cite{GL}.
\end{proof}

\subsection{Framing, symmetric colors and colored Jones specialization}
\label{sec:framing}
Changing the framing on a $j$-colored component of a link by $+1$ changes the colored HOMFLY-PT polynomial by a factor of $(-q)^{-j}a^{-j}q^{j^2}$, which is reflected on the generating functions for the skein module elements of these links by decreasing all entries for the exponent coefficient vectors $S$ and $A$ by one, while increasing all entries of the matrix $Q$ by one. In particular, after a suitable framing shift, the generating function in quiver form associated to a link may be assumed to have have a matrix $Q$ with non-negative integer values, in which case it is actually the adjacency matrix of a quiver. 

So far we have only been dealing with the \textit{anti-symmetric} HOMFLY-PT polynomials for links colored by one-column Young diagrams. However, we can use the following proposition to recover generating functions for \textit{symmetric} HOMFLY-PT polynomials with respect to one-row Young diagrams.

\begin{proposition}[{\cite[Proposition 4.4]{TVW}}] Let $L(\vec{\lambda})$ be a link with a coloring of its components described by the vector $\vec{\lambda}$ of Young diagrams. Then the unreduced colored HOMFLY-PT polynomial $\overline{P}$ satisfies
\[\overline{P}(L(\vec{\lambda^t}))=(-1)^c \overline{P}(L(\vec{\lambda^t}))|_{q\mapsto q^{-1}}\] 
where the superscript $t$ denotes transposition and $c$ is the total number of nodes in the Young diagrams in $\vec{\lambda}$.
\end{proposition}
\noindent This proposition implies that the one-row and one-column colored reduced HOMFLY-PT polynomials of knots can be obtained from each other by inverting $q$.  This operation changes the sign of all entries of the vector $S$ and of all entries of the matrix $Q$. Finally, one also needs to subtract one from every off-diagonal term in $Q$ to account for the asymmetry of $q$-Pochhammer symbols.

\begin{remark} The colored Jones polynomials can be recovered from the HOMFLY-PT polynomials with respect to colorings by one-row Young diagrams under the specialization $a\mapsto q^2$. Thus, Theorem~\ref{thm:links} implies that also the colored Jones polynomials of rational links admit generating functions in quiver form.
\end{remark}

\section{Examples}
\label{sec:examples}

\subsection{Torus knots}
We now use the explicit transformations from Section~\ref{sec:prooflemmas} to compute the quivers associated to torus knots $T_{2,2n+1}$. To this end we start with the generating function for the (rescaled) skein module element of the trivial upward oriented tangle:
\[  \scalemath{0.75}{ UP, \left( \begin{array}{c|c|c} 0 & 0 & 0 
\end{array} \right)},
 \scalemath{0.75}{
 \left( 
\begin{array}{c} 0 
\end{array} \right)
\quad\quad
\begin{tikzpicture} [scale=.5,anchorbase]
\draw[thin,opacity=.5] (-1,0) to (0,0);
\draw[thin,opacity=.5] (0,-.5) to (0,0);
\node at (0,0) {$\circ$};
\node at (-1.75,0) { $a^0$};
\node at (0,-1) {$q^{0}$};
\end{tikzpicture} 
} 
\]  
Here it is understood that the single summation index is inactive. The diagram on the right-hand side shows this summation index as a hollow dot, drawn in the plane at coordinates given by its $q$- and $a$-exponents when the color is specialized to $1$. Now we apply $TT$ and get:
\begin{equation}\label{eqn:ex1}  \scalemath{0.75}{ UP, \left( \begin{array}{c|c|c} 1 & -1 & 0 \\ \hline\hline 0 & -2 & 0 
\end{array} \right)},
 \scalemath{0.75}{
 \left( 
\begin{array}{c c} 0 & 0\\ 0 & 0 
\end{array} \right)
\quad\quad
\begin{tikzpicture} [scale=.5,anchorbase]
\draw[thin,opacity=.5] (-3,0) to (0,0);
\draw[thin,opacity=.5] (-2,-.5) to (-2,0);
\draw[thin,opacity=.5] (-1,-.5) to (-1,0);
\draw[thin,opacity=.5] (0,-.5) to (0,0);
\node at (-1,0) {$\bullet$};
\node at (-2,0) {$\circ$};
\node at (-3.75,0) { $a^0$};
\node at (-2,-1) { $q^{-2}$};
\node at (0,-1) {$q^{0}$};
\end{tikzpicture} 
} 
\end{equation}
Here we have two summation indices, one of which is active and drawn as a solid dot in the plane.
If $n=1$, we can now add the last crossing and perform the closure by applying the transformation from Lemma~\ref{lemmaclosure} and we get the following generating function data for the trefoil:
\[  \scalemath{0.75}{ \left( \begin{array}{c|c} -1 & -1 \\ \hline -3 & -1 \\ \hline  0 & 1 
\end{array} \right)},
 \scalemath{0.75}{
 \left( 
\begin{array}{c|c|c} 3 & 1 & 1\\ \hline 1 & 1 & 0 \\ \hline 1 & 0 & 0 
\end{array} \right) 
\quad\quad
\begin{tikzpicture} [scale=.5,anchorbase]
\draw[thin,opacity=.5] (-3,1) to (2,1);
\draw[thin,opacity=.5] (-3,0) to (2,0);
\draw[thin,opacity=.5] (-3,-1) to (2,-1); 
\draw[thin,opacity=.5] (-2,-1.5) to (-2,1);
\draw[thin,opacity=.5] (-1,-1.5) to (-1,1);
\draw[thin,opacity=.5] (0,-1.5) to (0,1);
\draw[thin,opacity=.5] (1,-1.5) to (1,1);
\draw[thin,opacity=.5] (2,-1.5) to (2,1);
\node at (-2,-1) {$\bullet$};
\node at (0,1) {$\bullet$};
\node at (2,-1) {$\bullet$};
\node at (-3.75,-1) {$a^{-1}$};
\node at (-3.75,1) { $a^1$};
\node at (-2,-2) { $q^{-2}$};
\node at (0,-2) {$q^0$};
\node at (2,-2) {$q^2$};
\draw[thick] (2.2,-1) circle (.2);
\draw[thick] (2.3,-1) circle (.3);
\draw[thick] (2.4,-1) circle (.4);
\draw[thick] (-2.4,-1) circle (.4);
\draw[thick] (-2,-1)to(2,-1);
\draw[thick] (0,1)to (2,-1);
\end{tikzpicture} 
}
\] 
The diagram on the right-hand side shows the quiver corresponding to this generating function. The vertices correspond to summation indices and they are drawn in coordinates determined by the corresponding $q$- and $a$-exponents when the color is specialized to one. 
 
In order to compare with the result in \cite{KRSSlong}, we have to adjust the framing by $3$, which shifts every entry in the exponent coefficient vectors for $-q$ and $a$ up by $3$ and every entry of the matrix $Q$ down by $3$.  
\[  \scalemath{0.75}{ \left( \begin{array}{c|c} 2 & 2 \\ \hline 0 & 2 \\ \hline  3 & 4 
\end{array} \right)},
 \scalemath{0.75}{
 \left( 
\begin{array}{c|c|c} 0 & -2 & -2\\ \hline -2 & -2 & -3 \\ \hline -2 & -3 & -3 
\end{array} \right)
} \]
Finally, we invert $q$ to proceed from the generating function for anti-symmetrically colored HOMFLY-PT polynomials to symmetrically colored HOMFLY-PT polynomials. For the trefoil the result is 
\[  \scalemath{0.75}{ \left( \begin{array}{c|c} -2 & 2 \\ \hline 0 & 2 \\ \hline  -3 & 4 
\end{array} \right)},
 \scalemath{0.75}{
 \left( 
\begin{array}{c|c|c} 0 & 1 & 1\\ \hline 1 & 2 & 2 \\ \hline 1 & 2 & 3 
\end{array} \right)
\quad\quad
\begin{tikzpicture} [scale=.5,anchorbase]
\draw[thin,opacity=.5] (-3,1) to (2,1);
\draw[thin,opacity=.5] (-3,0) to (2,0);
\draw[thin,opacity=.5] (-3,-1) to (2,-1); 
\draw[thin,opacity=.5] (-2,-1.5) to (-2,1);
\draw[thin,opacity=.5] (-1,-1.5) to (-1,1);
\draw[thin,opacity=.5] (0,-1.5) to (0,1);
\draw[thin,opacity=.5] (1,-1.5) to (1,1);
\draw[thin,opacity=.5] (2,-1.5) to (2,1);
\node at (-2,-1) {$\bullet$};
\node at (0,1) {$\bullet$};
\node at (2,-1) {$\bullet$};
\node at (-3.75,-1) {$a^2$};
\node at (-3.75,1) { $a^4$};
\node at (-2,-2) { $q^{-2}$};
\node at (0,-2) {$q^0$};
\node at (2,-2) {$q^2$};
\draw[thick] (0,1.2) circle (.2);
\draw[thick] (0,1.3) circle (.3);
\draw[thick] (0,1.4) circle (.4);
\draw[thick] (2.4,-1) circle (.4);
\draw[thick] (2.25,-1) circle (.25);
\draw[thick] (-2,-1)to(2,-1);
\draw[thick] (0,1)to (-2,-1);
\draw[thick] (0,1.075)to (2,-.925);
\draw[thick] (0,.925)to (2,-1.075);
\end{tikzpicture} 
} \]
which matches perfectly with Section 5.2 \cite{KRSSlong}.
 
For higher order torus knots we continue from \eqref{eqn:ex1} with further applications of the $TT$ operation before closing. Up to a permutation of summation indices and a framing change, the results again agrees with the color-transpositions of the quivers computed in Section 5.2 \cite{KRSSlong}.

\subsection{Knot \texorpdfstring{$\mathbf{7_3}$}{7_3}}
Quivers for all knots with at most six crossings, as well as all $(2,2n+1)$ torus  knots and twist knots have been computed in \cite{KRSSlong}. Representative for all rational knots outside of this set, we describe a quiver for the rational knot $7_3=K_{13/3}$ with the continued fraction expansion $[1,2,4]$. Here we start with the operation $RT$, then apply $TR$ and then $TT$, before executing the last top crossing and the closure. We show these steps, but omit to draw the arrows in the corresponding quivers. Applying $RT$ to the trivial upward oriented tangle produces:
\[  \scalemath{0.75}{ OP, \left( \begin{array}{c|c|c} 
0 & 0 & 0\\ \hline \hline 1& -2 & -1 
\end{array} \right)},
 \scalemath{0.75}{
 \left( 
\begin{array}{c|c} 
0 & 1\\ \hline 1 & 1 
\end{array} \right)
\quad\quad
\begin{tikzpicture} [scale=.5,anchorbase]
\draw[thin,opacity=.5] (-2,0) to (0,0);
\draw[thin,opacity=.5] (-2,-1) to (0,-1);
\draw[thin,opacity=.5] (0,-1.5) to (0,0);
\draw[thin,opacity=.5] (-1,-1.5) to (-1,0);
\node at (0,0) {$\bullet$};
\node at (-1,-1) {$\circ$};
\node at (-2.75,0) { $a^0$};
\node at (0,-2) {$q^{0}$};
\end{tikzpicture} 
} 
\]  
Then we apply $TR$ and get:
\[  \scalemath{0.75}{ UP, \left( \begin{array}{c|c|c} 
1&-1 & 0\\ \hline 1& -3 & -2\\ \hline \hline 0&-2 & 0\\ \hline 0&-4 & -2\\ \hline 0&-3 & -2
\end{array} \right)},
 \scalemath{0.75}{
 \left( 
\begin{array}{c|c|c|c|c} 
0 & 1 & 0 & 1 & 2\\ \hline 1 & 2 & 2 & 2 & 3\\ \hline 0 & 2 & 0 & 1 & 1\\ \hline 1 & 2 & 1 & 2 & 2\\ \hline 2 & 3 & 1 & 2 & 3
\end{array} \right)
\quad\quad
\begin{tikzpicture} [scale=.5,anchorbase]
\draw[thin,opacity=.5] (-3,0) to (0,0);
\draw[thin,opacity=.5] (-3,-1) to (0,-1);
\draw[thin,opacity=.5] (-3,-2) to (0,-2);
\draw[thin,opacity=.5] (0,-2.5) to (0,0);
\draw[thin,opacity=.5] (-1,-2.5) to (-1,0);
\draw[thin,opacity=.5] (-2,-2.5) to (-2,0);
\node at (-1,0) {$\bullet$};
\node at (-2,0) {$\circ$};
\node at (-1,-2) {$\bullet$};
\node at (-2,-2) {$\circ$};
\node at (0,-2) {$\circ$};
\node at (-3.75,0) { $a^0$};
\node at (-3.75,-2) { $a^{-2}$};
\node at (0,-3) {$q^{0}$};
\node at (-2,-3) {$q^{-2}$};
\end{tikzpicture} 
} 
\]  
Next, we apply $TT$ and get:
\[  \scalemath{0.75}{ UP, \left( \begin{array}{c|c|c} 
1 &-1 & 0\\ \hline 1& -3 & -2\\ \hline 1& -3 & 0\\ \hline 1& -5 & -2\\ \hline 1 & -4 & -2\\ \hline \hline 
0 & -4 & 0\\ \hline 0 & -6 & -2\\ \hline 0 & -5 & -2\end{array} \right)},
 \scalemath{0.75}{
 \left( 
\begin{array}{c|c|c|c|c|c|c|c} 
2 & 3 & 0 & 1 & 2 & 0 & 1 & 2\\ \hline 3 & 4 & 2 & 2 & 3 & 2 & 2 & 3\\ \hline 0 & 2 & 0 & 1 & 1 & 0 & 1 & 1\\ \hline 1 & 2 & 1 & 2 & 2 & 2 & 2 & 2\\ \hline 2 & 3 & 1 & 2 & 3 & 2 & 3 & 3\\ \hline 0 & 2 & 0 & 2 & 2 & 0 & 1 & 1\\ \hline 1 & 2 & 1 & 2 & 3 & 1 & 2 & 2\\ \hline 2 & 3 & 1 & 2 & 3 & 1 & 2 & 3
\end{array} \right)
\quad\quad
\begin{tikzpicture} [scale=.5,anchorbase]
\draw[thin,opacity=.5] (-5,0) to (1,0);
\draw[thin,opacity=.5] (-5,-1) to (1,-1);
\draw[thin,opacity=.5] (-5,-2) to (1,-2);
\draw[thin,opacity=.5] (1,-2.5) to (1,0);
\draw[thin,opacity=.5] (0,-2.5) to (0,0);
\draw[thin,opacity=.5] (-1,-2.5) to (-1,0);
\draw[thin,opacity=.5] (-2,-2.5) to (-2,0);
\draw[thin,opacity=.5] (-3,-2.5) to (-3,0);
\draw[thin,opacity=.5] (-4,-2.5) to (-4,0);
\node at (1,0) {$\bullet$};
\node at (1,-2) {$\bullet$};
\node at (-3,0) {$\bullet$};
\node at (-3,-2) {$\bullet$};
\node at (-1,-2) {$\bullet$};
\node at (-4,0) {$\circ$};
\node at (-4,-2) {$\circ$};
\node at (-2,-2) {$\circ$};
\node at (-5.75,0) { $a^0$};
\node at (-5.75,-2) { $a^{-2}$};
\node at (0,-3) {$q^{0}$};
\node at (-2,-3) {$q^{-2}$};
\node at (-4,-3) {$q^{-4}$};
\end{tikzpicture} 
} 
\]  
Finally, we display the result of the closure operation, followed by inverting $q$ and a framing adjustment:

\[\scalemath{0.75}{ \left( \begin{array}{c|c} 
-6 & 6\\ \hline -4 & 4\\ \hline -4 & 6\\ \hline -2 & 4\\ \hline -3 & 4\\ \hline -2 & 6\\ \hline 0 & 4\\ \hline -1 & 4\\ \hline -7 & 8\\ \hline -5 & 6\\ \hline -5 & 8\\ \hline -3 & 6\\ \hline -4 & 6
\end{array} \right)},
 \scalemath{0.75}{
 \left( 
\begin{array}{c|c|c|c|c|c|c|c|c|c|c|c|c} 
2 & 0 & 3 & 2 & 1 & 5 & 4 & 3 & 3 & 2 & 5 & 4 & 3\\ \hline 0 & 0 & 1 & 1 & 0 & 3 & 3 & 2 & 1 & 1 & 3 & 3 & 2\\ \hline 3 & 1 & 4 & 2 & 2 & 5 & 4 & 4 & 4 & 2 & 5 & 4 & 4\\ \hline 2 & 1 & 2 & 2 & 1 & 3 & 3 & 3 & 3 & 2 & 3 & 3 & 3\\ \hline 1 & 0 & 2 & 1 & 1 & 3 & 2 & 2 & 2 & 1 & 3 & 2 & 2\\ \hline 5 & 3 & 5 & 3 & 3 & 6 & 4 & 4 & 6 & 4 & 6 & 4 & 4\\ \hline 4 & 3 & 4 & 3 & 2 & 4 & 4 & 3 & 5 & 4 & 5 & 4 & 3\\ \hline 3 & 2 & 4 & 3 & 2 & 4 & 3 & 3 & 4 & 3 & 5 & 4 & 3\\ \hline 3 & 1 & 4 & 3 & 2 & 6 & 5 & 4 & 5 & 3 & 6 & 5 & 4\\ \hline 2 & 1 & 2 & 2 & 1 & 4 & 4 & 3 & 3 & 3 & 4 & 4 & 3\\ \hline 5 & 3 & 5 & 3 & 3 & 6 & 5 & 5 & 6 & 4 & 7 & 5 & 5\\ \hline 4 & 3 & 4 & 3 & 2 & 4 & 4 & 4 & 5 & 4 & 5 & 5 & 4\\ \hline 3 & 2 & 4 & 3 & 2 & 4 & 3 & 3 & 4 & 3 & 5 & 4 & 4
\end{array} \right)
\quad\quad
\begin{tikzpicture} [scale=.5,anchorbase]
\draw[thin,opacity=.5] (-5,4) to (4,4);
\draw[thin,opacity=.5] (-5,5) to (4,5);
\draw[thin,opacity=.5] (-5,6) to (4,6);
\draw[thin,opacity=.5] (-5,7) to (4,7);
\draw[thin,opacity=.5] (-5,8) to (4,8);
\draw[thin,opacity=.5] (3.5,8) to (-1,8);
\draw[thin,opacity=.5] (-4,3.5) to (-4,8);
\draw[thin,opacity=.5] (-3,3.5) to (-3,8);
\draw[thin,opacity=.5] (-2,3.5) to (-2,8);
\draw[thin,opacity=.5] (-1,3.5) to (-1,8);
\draw[thin,opacity=.5] (0,3.5) to (0,8);
\draw[thin,opacity=.5] (1,3.5) to (1,8);
\draw[thin,opacity=.5] (2,3.5) to (2,8);
\draw[thin,opacity=.5] (3,3.5) to (3,8);
\draw[thin,opacity=.5] (4,3.5) to (4,8);
\node at (-4,6) {$\bullet$};
\node at (-4,4) {$\bullet$};
\node at (0,6) {$\bullet$};
\node at (0,4) {$\bullet$};
\node at (-2,4) {$\bullet$};
\node at (4,6) {$\bullet$};
\node at (4,4) {$\bullet$};
\node at (2,4) {$\bullet$};
\node at (-2,8) {$\bullet$};
\node at (-2,6) {$\bullet$};
\node at (2,8) {$\bullet$};
\node at (2,6) {$\bullet$};
\node at (0,6) {$\bullet$};
\node at (-5.75,4) { $a^4$};
\node at (-5.75,6) { $a^{6}$};
\node at (-5.75,8) { $a^{8}$};
\node at (0,3) {$q^{0}$};
\node at (-2,3) {$q^{-2}$};
\node at (-4,3) {$q^{-4}$};
\node at (2,3) {$q^{2}$};
\node at (4,3) {$q^{4}$};
\end{tikzpicture} 
}\]

\subsection{More examples}
Using a SageMath implementation of the algorithm described in Section~\ref{sec:knots}, we have computed the generating function data for all 362 rational knots with up to twelve crossings, which in total takes less than a minute on a standard PC. This data as well as the source code is available upon request. 
\medskip

\noindent\textbf{Funding:} The work of M.~S. was supported by the European Research Council [Starting Grant no. 335739] ''Quantum fields and knot homologies" funded by the European Research Council under the European Union's Seventh Framework Programme, and by the Ministry of Education, Science, and Technological Development of the Republic of Serbia [project no. 174012]. The work of P.~W. was supported by the Leverhulme Trust [Research Grant RP2013-K-017] and the Australian Research Council Discovery Projects ``Braid groups and higher representation theory'' and ``Low dimensional categories'' [DP140103821, DP160103479].
\medskip

\noindent\textbf{Acknowledgments:} We are grateful to the anonymous referees whose suggestions have significantly improved the presentation of the paper. We would like to thank Jake Rasmussen and Piotr Su{\l}kowski for valuable discussions. Part of this work was done at the Max Planck Institute for Mathematics, Bonn, Germany, and at the Isaac Newton Institute for Mathematical Sciences, Cambridge, UK, during the programme ``Homology theories in low dimensional topology'', which was supported by the UK Engineering and Physical Sciences Research Council [Grant Number EP/K032208/1]. We thank both institutions for their hospitality.

This is a pre-copyedited, author-produced version of an article accepted for publication in International Mathematical Research Notices following peer review. The version of record is:
\\
Rational Links and DT Invariants of Quivers,
Marko Sto$\check{\text{s}}$i$\acute{\text{c}}$ and Paul Wedrich,
International Mathematics Research Notices, rny289, \url{https://doi.org/10.1093/imrn/rny289}.

\end{document}